\documentclass[11pt]{article}
\usepackage[top = 1.3in, bottom= 1.3in, left = 1.3in, right=1.3in]{geometry}

\usepackage{paper2-macros}
\title{Ideal triangulations and geometric transitions}
\author{Jeffrey Danciger\footnote{J.D. is partially supported by the National Science Foundation, grant number DMS 1103939}\\
  Department of Mathematics,\\
  University of Texas - Austin,\\
  1 University Station C1200,\\
  Austin, TX 78712-0257\\
  \texttt{jdanciger@math.utexas.edu}}
\date{\today}

\begin{document}

\maketitle

\begin{abstract}
Thurston introduced a technique for finding and deforming three-dimensional hyperbolic structures by gluing together ideal tetrahedra. 
We generalize this technique to study families of geometric structures that transition from hyperbolic to anti de Sitter (AdS) geometry.
Our approach involves solving Thurston's gluing equations over several different shape parameter algebras. 
In the case of a punctured torus bundle with Anosov monodromy, we identify two components of real solutions for which there are always nearby positively oriented solutions over both the complex and pseudo-complex numbers. These complex and pseudo-complex solutions define hyperbolic and AdS structures that, after coordinate change in the projective model, may be arranged into one continuous family of real projective structures. We also study the rigidity properties of certain AdS structures with tachyon singularities.

\end{abstract}

\section{Introduction}

In his notes, Thurston showed how to construct and deform hyperbolic structures by glueing together hyperbolic ideal tetrahedra (\cite{Thurston} or see \cite{Neumann-85}). This technique is often used to study the hyperbolic Dehn surgery space for a three-dimensional manifold $M$ with a union of tori as boundary. Supposing $M$ is equipped with a topological ideal triangulation, the basic idea is to assign a hyperbolic structure, or \emph{shape}, to each tetrahedron, so that the tetrahedra fit together consistently to give a smooth hyperbolic structure on~$M$. 
The shape of a tetrahedron is described by a complex variable, and the consistency conditions are algebraic equations, known as Thurston's equations or the \emph{gluing equations}. 
A solution defines a hyperbolic structure on $M$ if all of the tetrahedra are oriented consistently. However, there are many solutions for which some tetrahedra are inverted or collapsed. In many cases, there exist solutions for which all of the tetrahedra are collapsed onto a plane; these are the solutions for which all shapes are real numbers.
This article explores the geometric significance of these real solutions, each of which defines a transversely hyperbolic foliation~$\mathcal F$, meaning a one-dimensional foliation of~$M$ with a hyperbolic structure transverse to the leaves (see~\cite{Thurston} Definition~4.8.2). In particular, we study the question of \emph{regeneration}, first studied by Hodgson in \cite{Hodgson-86}: Is $\mathcal F$ the limit of a path of collapsing hyperbolic structures~$\mathscr H_t$? If the answer is yes, we say that $\mathcal F$ regenerates to hyperbolic structures. In this triangulated setting, we ask further for the hyperbolic structures $\mathscr H_t$ to be constructed from nearly collapsed positively oriented ideal tetrahedra which limit to the collapsed tetrahedra that define~$\mathcal F$.

Anti de Sitter ($\AdS$) geometry is a Lorentzian analogue of hyperbolic geometry; it is the homogeneous Lorentzian geometry of constant negative curvature. 
Recently, in \cite{Danciger-13-1}, the author demonstrated that transversely hyperbolic foliations may also arise as limits of collapsing $\AdS$ structures. 
Therefore, we also study the regeneration question in the AdS context: Is a given transversely hyperbolic foliation $\mathcal F$ the limit of collapsing anti de Sitter structures $\mathscr A_t$? We use Thurston's technique to study the regeneration question in both the hyperbolic and $\AdS$ context using the same algebraic framework. For we show that (incomplete) triangulated $\AdS$ structures on~$M$ may be constructed by solving the gluing equations over the \emph{pseudo-complex} numbers (Section~\ref{AdS-tetrahedra}). The pseudo-complex numbers, generated by a non-real element $\tau$ such that $\tau^2 = +1$, has a Lorentzian similarity structure analogous to the Euclidean similarity structure on the complex numbers. We note that a pseudo-complex number $z$ determines an ideal tetrahedron in $\AdS$ if and only if $z$, $1/(1-z)$, and $(z-1)/z$ are defined and spacelike.

Our main result is that when $M$ is a punctured torus bundle with Anosov monodromy, many transversely hyperbolic foliations do regenerate to both hyperbolic and $\AdS$ structures. Further, the hyperbolic and $\AdS$ structures regenerated from a given transversely hyperbolic foliation fit naturally into a continuous path of geometric structures, giving new examples of the geometric transition phenomenon described in \cite{Danciger-13-1}. Finally, we apply the main result to study the space of $\AdS$ structures on $M$ whose geodesic completion has a special singularity called a tachyon.

\subsection{Regeneration results} \label{intro-regen}
Throughout, let $M$ be a three-dimensional manifold with a union of tori as boundary and suppose $\mathcal{T} = \{\mathcal T_1,\ldots,\mathcal T_N\}$ is a topological ideal triangulation of $M$. The space $\mathscr D_\RR$ of real solutions to Thurston's equations will be called the \emph{real deformation variety}. Each point of $\mathscr D_\RR$ determines a transversely hyperbolic foliation on $M$ built by projecting the simplices of~$\mathcal T$ onto ideal quadrilaterals in $\HH^2$; such a structure will be called a transversely hyperbolic foliation on $(M, \mathcal T)$.
In Section~\ref{s:tet-regen}, we observe simple conditions (Proposition~\ref{tet-regen}) on the real deformation variety $\mathscr D_\RR$ guaranteeing regeneration to hyperbolic and $\AdS$ structures. Specifically, if $(z_j) \in \RR^N$ is a solution to Thurston's equations defining a transversely hyperbolic foliation $\mathcal F$, then $\mathcal F$ regenerates to hyperbolic and $\AdS$ structures~if 
\begin{enumerate}
\itemsep0em
\item[(1)]  $\mathscr D_\RR$ is smooth at $(z_j)$ and 
\item[(2)]  there exists a tangent vector $(v_j) \in \RR^N$ to $\mathscr D_\RR$ at $(z_j)$ such that $v_j > 0$ for all $j$.
\end{enumerate}
A tangent vector $(v_j)$ as above is called a \emph{positive tangent vector}.

 In the case that $M$ is a punctured torus bundle with Anosov monodromy, we identify two canonical connected components of the real deformation variety for which every point has the best possible properties. The following theorem, discussed in Section~\ref{torus-bundles}, is the fundamental result of the paper, on which the main theorems (\ref{thm:regen}, \ref{triangulated-transition} and~\ref{thm:tachyons} below) rest.
\begin{Theorem}\label{punctured-torus}
Let $M^3$ be a punctured torus bundle over the circle with Anosov monodromy and let $\mathcal T$ be the monodromy triangulation of $M$. Let $\mathscr D_\RR$ be the deformation variety of real solutions to Thurston's equations for $(M, \mathcal{T})$. Then, there are two canonical smooth, one dimensional, connected components $\Vpl$ of $\mathscr D_\RR$ with positive tangent vectors at every point. Each component of $\Vpl$ is parameterized by the (signed) length of the puncture curve.
\end{Theorem}
\noindent Making use of Gu\'{e}ritaud's description of Thurston's equations for the monodromy triangulation (see \cite{Gueritaud-06}), we prove Theorem~\ref{punctured-torus} by direct analysis of the equations.
Via the discussion above, Theorem~\ref{punctured-torus} implies\marginnote{change from cor to thm thoughout}:
\begin{Theorem}
\label{thm:regen}
Let $\mathcal F$ be a transversely hyperbolic foliation on $(M, \mathcal T)$ determined by a point of $\Vpl$. Then there exist hyperbolic structures $\mathscr H_t$, and $\AdS$ structures $\mathscr A_t$ on $(M, \mathcal T)$, defined for $t > 0$, such that $\mathscr H_t$ and $\mathscr A_t$ collapse to $\mathcal F$ as $t \to 0$.
\end{Theorem}

Let us briefly describe the relationship between the two components of $\Vpl$, denoted by~$\Vpl^+$ and~$\Vpl^-$. We construct two special solutions $(z_j^+)$ and $(z_j^-)$ which serve as base-points for $\Vpl^+$ and $\Vpl^-$ respectively. These two solutions come from the Sol geometry of the torus bundle obtained by filling in the puncture. Following work of Huesener--Porti--Su\'arez \cite{Huesener-01}, we show that though $(z_j^+)$ and $(z_j^-)$ determine distinct transversely hyperbolic foliations $\mathcal F^+$ and $\mathcal F^-$, the transversely hyperbolic foliations determined by $\Vpl^+ \setminus \{ (z^+_j) \}$ are equivalent to those determined by  $\Vpl^- \setminus \{ (z^-_j) \}$ (though the triangulations spin in opposite directions). Thus the corresponding component of the space of transversely hyperbolic foliations on $M$ (forgetting the triangulation) is a  ``line with two origins" (the two origins are $\mathcal F_+$ and $\mathcal F_-$); in particular the deformation space is not Hausdorff.

\subsection{Geometric transitions}

In \cite{Danciger-13-1}, the author defined a notion of geometric transition from hyperbolic to AdS structures in the context of real projective geometry. Specifically, a transition is a path $\mathscr P_t$ of real projective structures on a manifold $M$ such that $\mathscr P_t$ is conjugate to a hyperbolic structure if $t > 0$, or to an $\AdS$ structure if $t < 0$. The author then constructed the first examples of this hyperbolic-AdS transition in the case that the manifold $M$ is the unit tangent bundle of the $(2,m,m)$ triangle orbifold (note that the structures have cone singularities). 
In Section~\ref{tet-transitions} we show how to build examples of the hyperbolic-$\AdS$ transition in the triangulated setting when the conditions~(1) and~(2) above\marginnote{these are conditions on a particular real solution, which is not referenced in this sentence...} are satisfied.
In particular, the geometric structures produced by Theorem~\ref{thm:regen} may be organized into a geometric transition:

\begin{Theorem}\label{triangulated-transition}
Let $M^3$ be an Anosov punctured torus bundle, and let $\Vpl$ be as in Theorem~\ref{punctured-torus}. Given a transversely hyperbolic foliation $\mathcal F$ determined by a point of $\Vpl$, let $\mathscr H_t$ and $\mathscr A_t$ be the hyperbolic and $\AdS$ structures that collapse to $\mathcal F$, given by Theorem~\ref{thm:regen}. Then, there exists a continuous path $\mathscr P_t$ of real projective structures on $M$ such that 
\begin{itemize}
\itemsep0em
\item $\mathscr P_t$ is conjugate to the hyperbolic structure $\mathscr H_t$ if $t > 0$.
\item $\mathscr P_t$ is conjugate to the $\AdS$ structure $\mathscr A_{|t|}$ if $t < 0$.
\end{itemize}
For each $t$ (including $t=0$), the monodromy triangulation $\mathcal T$ is realized in $\mathscr P_t$ by positively oriented tetrahedra.
\end{Theorem}

For $t > 0$ (resp. $t < 0$), the projective structures $\mathscr P_t$ should be thought of as rescaled versions of the $\mathscr H_t$ (resp. the $\mathscr A_{|t|}$); they are obtained by applying a linear transformation of $\RP^3$ that stretches $\HH^3$ (and $\AdS^3$) transverse to the collapsing direction. The projective structure $\mathscr P_0$ obtained in the limit is a \emph{half-pipe} (HP) structure. Half-pipe geometry, defined in \cite{Danciger-13-1}, is the natural transitional geometry bridging the gap between hyperbolic and $\AdS$ geometry. An $\HP$ structure encodes information about the collapsed hyperbolic structure $\mathcal F$ as well as first order information in the direction of collapse.
The proof of the theorem leads naturally to the construction of $\HP$ structures out of $\HP$ ideal tetrahedra. The algebra describing the shapes of $\HP$ tetrahedra is the tangent bundle of~$\RR$, thought of as complex (or pseudo-complex) numbers with infinitesimal imaginary part.

\subsection{$\AdS$ structures with tachyons}
A \emph{tachyon} singularity (\cite{Barbot-09} or see \cite{Danciger-13-1}), in $\AdS$ geometry is a singularity along a space-like geodesic $\mathcal L$ such that the holonomy of a meridian encircling $\mathcal L$ is a Lorentz boost that point-wise fixes $\mathcal L$. The (signed) magnitude of the boost is called the \emph{mass}. Tachyons are the Lorentzian analogue of cone singularities in hyperbolic geometry, with the tachyon mass playing the role of the cone angle. Paths of geometric structures that transition from hyperbolic cone-manifolds to $\AdS$ manifolds with tachyons are constructed in \cite{Danciger-13-1}.
In Section~\ref{AdS-tetrahedra}, we describe how to construct $\AdS$ structures with tachyons by solving the gluing equations with an added condition governing the geometry at the boundary.  

Hodgson-Kerckhoff \cite{Hodgson-98} showed that, under assumptions about the cone angle, compact hyperbolic cone-manifolds are parameterized by the cone angles. In particular, hyperbolic cone-manifolds are rigid relative to the cone angles. It is natural to ask whether a similar rigidity phenomenon occurs in the $\AdS$ setting (see \cite{BBDGGKKSZ-12}). We give some evidence to the affirmative in Section~\ref{sec:tachyon}:

\begin{Theorem}
\label{thm:tachyons}
Let $N$ be a torus bundle over the circle with Anosov monodromy, and let $\Sigma$ be the simple closed curve going once around the base circle direction. Then the space of $\AdS$ structures on $N$ with a tachyon along $\Sigma$ contains a smooth one-dimensional connected component, parameterized by the tachyon mass. The mass can take any value in $(-\infty,0)$.
\end{Theorem}

 There is no analogue of Mostow rigidity in the $\AdS$ setting. Indeed, compact $\AdS$ manifolds without singularities are very flexible. However, Theorem~\ref{thm:tachyons} suggests that the presence of a tachyon singularity may restrict the geometry significantly, so that the deformation theory mimics that of hyperbolic cone-manifolds.

\subsection{Some related literature}

Following Thurston's pioneering work, ideal triangulations have proven an important theoretical and experimental tool in hyperbolic geometry. The software packages SnapPea \cite{SnapPea}, Snap \cite{Snap}, and SnapPy \cite{SnapPy}, frequently used for experiment in hyperbolic three-manifolds, find hyperbolic structures via ideal triangulations. 
The volume maximization program of Casson-Rivin~\cite{Rivin-94}, the recent variational formulation of the Poincar\'e conjecture by Luo~\cite{Luo-10}, and many other articles (e.g. \cite{Lackenby-00, Hodgson-12}) feature the more general notion of angled ideal triangulation.  Further, the use of ideal triangulations has spread beyond hyperbolic geometry, including representation theory and quantum topology (see e.g. \cite{Zickert-09, Dimofte-13}). However, the present work gives, to the author's knowledge, the first application of ideal triangulations in the context of $\AdS$ geometry. Our results add to the growing list of parallels in the studies of hyperbolic and $\AdS$ geometry, beginning with Mess's classification of globally hyperbolic $\AdS$ space-times \cite{Mess-07, Andersson-07} and its remarkable similarity to the Simultaneous Uniformization theorem of Bers \cite{Bers-60}. 
Stemming from Mess's work, results and questions in hyperbolic and $\AdS$ geometry (see \cite{BBDGGKKSZ-12}) have begun to appear in tandem, suggesting the existence of a deeper link. 
The geometric transition construction of \cite{Danciger-13-1} and its triangulated counterpart described here give a concrete connection strengthening the analogy between the two geometries.

\bigskip

{\noindent \bf Organization of the paper.} 
The paper naturally divides into two parts. The first part, consisting of Sections~\ref{models}--\ref{tet-transitions}, describes a general algebraic framework for building triangulated geometric structures, including $\HH^3$ and $\AdS$ structures, as well as structures that transition between the two geometries. Section~\ref{models} constructs the relevant homogeneous spaces, and sets the notation used throughout Sections~\ref{general-tets}--\ref{tet-transitions}. Section~\ref{general-tets} gives a general construction of ideal tetrahedra and describes Thurston's equations in this setting. Section~\ref{triangulated-structures} applies the general construction to the cases of hyperbolic, $\AdS$, and $\HP$ structures, as well as transversely hyperbolic foliations. Section~\ref{s:tet-regen} gives the proof for the regeneration conditions described above in \ref{intro-regen}. Section~\ref{tet-transitions} describes how geometric structures transitioning from hyperbolic to $\AdS$ may be built from ideal tetrahedra, and proves that Theorem~\ref{triangulated-transition} follows from Theorem~\ref{punctured-torus}.

The second part of the paper focuses on punctured torus bundles. In Section~\ref{torus-bundles}, we give a description of the monodromy triangulation and then carefully study the real deformation variety and prove Theorem~\ref{punctured-torus}. Section~\ref{torus-bundles} may be read independently from the first part of the paper. Section~\ref{sec:tachyon} gives the proof of Theorem~\ref{thm:tachyons}.

The paper develops much of the needed hyperbolic and $\AdS$ geometry directly, and for this reason there is no section dedicated to preliminaries. The reader may wish to consult \cite{Thurston, Ratcliffe-book} for basics on hyperbolic geometry, or \cite{Benedetti-09, Danciger-11, Goldman-13} for basics on $\AdS$ geometry.

\bigskip

{\noindent \bf Acknowledgements.} I am very grateful to Steven Kerckhoff for advising this work during my doctoral studies at Stanford University. In addition, many helpful discussions occurred at the \emph{Workshop on Geometry, Topology, and Dynamics of Character Varieties} at the National University of Singapore in July 2010. In particular, discussions with Francesco Bonsante, Craig Hodgson, Feng Luo, and Jean-Marc Schlenker were very helpful in developing the ideas presented here.

%
%

\section{The algebra $\mathcal B$ and the model space $\XX$}
\label{models}
In \cite{Danciger-13-1}, the author constructs a family of model geometries in projective space that transitions from hyperbolic geometry to $\AdS$ geometry, passing though \emph{half-pipe} geometry.
We review the dimension-three version of this construction here. 
Each model geometry $\XX$ is associated to a real two-dimensional commutative algebra $\BB$. In Section~\ref{general-tets}, we will see that (a subset of) $\BB$ describes the shapes of ideal tetrahedra in~$\XX$.

Let $\BB = \mathbb R + \mathbb R \kappa$ be the real two-dimensional, commutative algebra generated by a non-real element $\kappa$ with $\kappa^2 \in \RR$. As a vector space $\BB$ is spanned by $1$ and $\kappa$. There is a conjugation action: 
$ \overline{(a + b \kappa)} := a - b \kappa, $
which defines a square-norm
$$ |a + b\kappa|^2 := (a + b\kappa)\overline{(a + b\kappa)} = a^2 -b^2 \kappa^2 \ \in \ \mathbb R.$$
Note that $|\cdot|^2$ may not be positive definite. We refer to $a$ as the \emph{real part}  and $b$ as the \emph{imaginary part} of $a + b\kappa$. If $\kappa^2 = -1$, then our algebra $\mathcal B = \CC$ is just the complex numbers, and in this case we use the letter $i$ in place of $\kappa$, as usual. If $\kappa^2 = +1$, then $\mathcal B$ is the \emph{pseudo-complex numbers} and we use the letter $\tau$ in place of $\kappa$.  
Note that if $\kappa^2 < 0$, then $\mathcal B \cong \CC$, and if $\kappa^2 > 0$ then $\mathcal B \cong \Rtau$.
In the case $\kappa^2 = 0$, we use the letter $\sigma$ in place of $\kappa$. In this case $\mathcal B = \Rsigma$ is isomorphic to the tangent bundle of the real numbers.

Now consider the $2\times2$ matrices $M_2(\BB)$. 
Let $\Herm(2, \BB) = \{ A \in M_2(\BB) : A^* = A\}$ denote the $2\times 2$ Hermitian matrices, where $A^*$ is the conjugate transpose of $A$.  As a real vector space, $\Herm(2,\BB) \cong \mathbb R^4$. We define the following (real) inner product on $\Herm(2,\BB)$:
$$ \left\langle \bminimatrix{a}{z}{\bar{z}}{d}, \bminimatrix{e}{w}{\bar{w}}{h} \right\rangle = -\frac{1}{2}tr\left( \bminimatrix{a}{z}{\bar{z}}{d} \bminimatrix{h}{-w}{-\bar{w}}{e}\right).$$
We will use the coordinates on $\Herm(2,\BB)$ given by 
\begin{align} \label{coordinates-on-Herm}
X &= \bminimatrix{x_1+x_2}{x_3 - x_4 \kappa}{x_3 + x_4\kappa}{x_1 - x_2}.
\end{align} In these coordinates, we have that
$$\langle X, X\rangle = -\text{det}(X) = -x_1^2 + x_2^2 + x_3^2 -\kappa^2 x_4^2$$
and we see that the signature of the inner product is $(3,1)$ if $\kappa^2 < 0$, or $(2,2)$ if $\kappa^2 > 0$.

The coordinates above identify $\Herm(2,\BB)$ with $\RR^4$. Then, identify the real projective space $\RP^3$ with the non-zero elements of $\Herm(2,\BB)$, considered up to multiplication by a real number. We define the region $\XX$ inside $\RP^3$ as the negative lines with respect to $\langle \cdot, \cdot\rangle$:
$$\XX = \left\{ X \in \Herm(2,\BB) : \langle X, X\rangle < 0 \right\} / \RR^*$$
We define the group $\PGL^+(2,\mathcal B)$ to be the $2\times 2$ matrices $A$, with coefficients in $\mathcal B$, such that $|\det(A)|^2 > 0$, up to the equivalence $A \sim \lambda A$ for any $\lambda \in \mathcal B$.

The group $\PGL^+(2, \BB)$ acts on $\XX$ by orientation preserving projective linear transformations as follows. Given $A \in \PGL^+(2,\BB)$ and $X \in \XX$:
$$ A \cdot X := A X A^*.$$

Note that if $\BB = \CC$, then $\PGL^+(2,\BB) = \PSL(2,\CC)$ and $\XX$ identifies with the usual projective model for \emph{hyperbolic space} $\XX = \HH^3$. In this case, the action above is the usual action by orientation preserving isometries of $\HH^3$, and gives the familiar isomorphism $\PSL(2,\CC) \cong \PSO(3,1)$,

If $\BB = \Rtau$, with $\tau^2 = +1$, then $\XX$ identifies with the usual projective model for \emph{anti de Sitter space} $\XX = \AdS^3$. Anti de Sitter geometry is a Lorentzian analogue of hyperbolic geometry. The natural invariant metric on $\AdS^3$, inherited from $\langle \cdot, \cdot \rangle$, has signature $(2,1)$. As such tangent vectors are partitioned into three types: \emph{space-like}, meaning positive, $\emph{time-like}$, meaning negative, or \emph{light-like}, meaning null. In any given tangent space, the light-like vectors form a cone that partitions the time-like vectors into two components. Thus, locally there is a continuous map assigning the name \emph{future pointing} or \emph{past pointing} to time-like vectors. The space $\AdS^3$ is \emph{time-orientable}, meaning that the labeling of time-like vectors as future or past may be done consistently over the entire manifold. The action of $\PGL^+(2,\Rtau)$ on $\AdS^3$ is by isometries, thus giving an embedding $\PGL^+(2,\Rtau) \hookrightarrow \PSO(2,2)$. In fact, $\PGL^+(2,\Rtau)$ has two components, distinguished by whether or not the action on $\AdS^3$ preserves time-orientation, and the map is an isomorphism.

Lastly, we discuss the case $\BB = \Rsigma$, with $\sigma^2 = 0$. In this case, $\XX = \HP^3$ is the projective model for \emph{half-pipe geometry} $(\HP)$, defined in \cite{Danciger-13-1} for the purpose of describing a geometric transition going from hyperbolic to $\AdS$ structures. 
The algebra $\mathbb R + \mathbb R \sigma$ should be thought of as the tangent bundle of $\mathbb R$: Letting $x$ be the standard coordinate function on $\mathbb R$, we think of $a + b \sigma$ as a path based at $a$ with tangent $b \frac{\partial}{\partial x}$. This point of view is particularly appropriate in the context of geometric transitions. For, given collapsing hyperbolic structures with holonomy representations $\rho_t  : \pi_1 M \rightarrow \PSL(2,\mathbb C)$ converging to $\rho_0 :\pi_1 M \rightarrow \PSL(2,\mathbb R)$, there is a natural $\HP$ representation defined by $$\rho_{\HP}(\cdot) = \rho_0(\cdot) + \sigma \operatorname{Im}\rho'_0(\cdot)$$
where $\operatorname{Im} \rho'_0(\gamma) \in T_{\rho_0(\gamma)} \PSL(2,\mathbb R)$ is the imaginary part of the derivative of $\rho_t$.
 We show in \cite{Danciger-13-1} how, in certain cases, the collapsing hyperbolic structures corresponding to $\rho_t$ limit, as real projective structures, to a (non-collpased) half-pipe structure. In that case $\rho_{\HP}$ is the holonomy representation of this limiting structure.
A similar interpretation is possible in the context of collapsing $\AdS$ structures with holonomy representations $\rho_t : \pi_1 M \rightarrow \PSL(2,\mathbb R + \mathbb R \tau)$ whose imaginary parts are going to zero.

\begin{Remark}
In each case, the orientation reversing isometries are also described by $\PGL^+(2, \BB)$ acting by $X \mapsto A \overline{X} A^*$.
\end{Remark}

{\noindent \bf The ideal boundary.} The ideal boundary $\partial^\infty \XX$ is the boundary of the region $\XX$ in $\RP^3$. It is given by the null lines in $\Herm(2,\BB)$ with respect to $\langle \cdot, \cdot \rangle$. Thus $$\partial^\infty \XX = \left\{ X \in \Herm(2,\BB) : \det(X) = 0, X \neq 0\right\}/\RR^*$$ can be thought of as the $2\times 2$ Hermitian matrices of rank one. We now give a useful description of $\partial^\infty \XX$ that generalizes the identification $\partial^\infty \HH^3 = \CP^1$. 

Any rank one Hermitian matrix $X$ can be decomposed (up to $\pm$) as 
\begin{equation} \label{eqn:decomposition}
X= \pm v v^*
\end{equation} where $v \in \BB^2$ is a two-dimensional column vector with entries in $\BB$, unique up to multiplication by $\lambda \in \BB$ with $|\lambda|^2 = 1$ (and $v^*$ denotes the transpose conjugate). 
This gives the identification $$\partial ^\infty \mathbb X \cong \mathbb P^1 \BB = \left\{ v \in \BB^2 : v v^* \neq 0 \right\} / \sim$$ where $v \sim v \lambda \text{ for } \lambda \in \BB^{\times}.$  The action of $\PGL^+(2, \BB)$ on $\mathbb P^1 \BB$ by matrix multiplication extends the action of $\PGL^+(2,\BB)$ on $\XX$ described above.

We use the square-bracket notation $\btwovector{x}{y}$ to denote the equivalence class in $\mathbb P^1 \BB$ of $\begin{pmatrix}x \\ y \end{pmatrix} \in \BB^2$. Similarly, a $2\times 2$ square-bracket matrix $\bminimatrix{a}{b}{c}{d}$ denotes the equivalence class in $\PGL^+(2,\BB)$ of the matrix $\minimatrix{a}{b}{c}{d} \in \GL^+(2,\BB)$. Throughout, we will identify $\BB$ with its image under the injection $\BB \hookrightarrow \mathbb P^1 \BB$ given by $z \mapsto \btwovector{z}{1}$.

\begin{Remark}
In the case $\kappa^2 \geq 0$, the condition $v v* \neq 0$ in the definition of $\mathbb P^1 \BB$ is \emph{not} equivalent to the condition $v \neq 0$, because $\BB$ has zero divisors.
\end{Remark}

%
%
\section{General construction of ideal tetrahedra}

\label{general-tets}

The following construction, based on Thurston's description of ideal tetrahedra in $\HH^3$, takes place in the model space $\XX$ defined by a two-dimensional real commutative algebra $\BB$, from the previous section. The generality of this construction is its main advantage. Geometric interpretations for the cases of interest will be given in the following section.
 
 Let $\Delta^n$ be the standard $n$-simplex with vertices removed: $$\Delta^n = \left\{ (t_1, \ldots, t_{n+1}) \in \RR^{n+1} : 0 \leq t_i < 1, \sum_i t_i = 1\right\}.$$
 Then a choice of $n+1$ rank one Hermitian matrices $Z_1, \ldots, Z_{n+1}$ determines a map $i: \Delta^n \to \Herm(2,\BB) \to \RP^3$ by $$i: (t_1, \ldots, t_{n+1}) \mapsto [t_1 Z_1 + \cdots + t_{n+1} Z_{n+1}].$$
We say that $Z_1, \ldots, Z_{n+1}$ determines an \emph{ideal $n$-simplex} in $\XX$ if the image of $\Delta^n$ lies entirely in $\XX$. When $n = 2$, an ideal simplex is called an \emph{ideal triangle}, and when $n=3$, an ideal simplex is called an \emph{ideal tetrahedron}. Note that we allow $i(\Delta^n)$ to have dimension smaller than $n$; in this case the ideal simplex is called \emph{collapsed}.
 The following proposition is elementary.
\begin{Proposition}
$Z_1, \ldots, Z_{n+1}$ determine an ideal simplex if and only if
\begin{equation}
\langle Z_i, Z_j \rangle < 0 \text{ for all } i\neq j
\label{negativity}
\end{equation}
\end{Proposition}
It is apparent from the proposition that the ideal simplex $i:\Delta^n \to \XX$ determined by $Z_1, \ldots, Z_{n+1}$ depends only on the equivalence classes $[Z_1], \ldots, [Z_{n+1}]$, up to reparametrization of $\Delta^n$. By the discussion in the previous section, we identify the ideal vertices $[Z_i] \in \partial \XX$ with elements $z_i \in \mathbb P^1 \BB$, and refer to $i: \Delta^n \to \XX$ as the tetrahedron determined by $z_1, \ldots, z_{n+1}$.

The following generalizes the three-transitivity of the action of $\PSL(2,\CC)$ on $\CP^1$.
\begin{Proposition}
\label{prop:standard-position}
Let $Z_1, Z_2, Z_3$ be rank one Hermitian matrices that determine an ideal triangle. Let $z_1, z_2, z_3$ be the corresponding elements of $\mathbb P^1 \BB$. Then there exists a unique $A \in \PGL^+(2,\BB)$ placing the $z_i$ in \emph{standard position}:  $$A z_1 = \infty = \btwovector{1}{0}, A z_2 = 0 = \btwovector{0}{1}, A z_3 = 1 = \btwovector{1}{1}.$$
\end{Proposition}
\begin{proof}
By transitivity of the action of $\PGL^+(2,\BB)$ on $\mathbb P^1 \BB$, we may assume that $z_1 = \infty$. Further, assume that $Z_1 = \minimatrix{1}{0}{0}{0}$. Writing $Z_2 = \minimatrix{a}{b}{\bar b}{d}$, with $a,d \in \RR$, $b \in \BB$ and $a d = b \bar b$, we have that $\langle Z_1, Z_2 \rangle = -\frac{d}{2}.$ It follows that $d > 0$, and so $d \cdot Z_2 = \twovector{b}{d} \begin{pmatrix}\bar b & d \end{pmatrix}$, in other words $z_2 = \frac{b}{d} \in \BB$. Then, the element $\bminimatrix{1}{-z_2}{0}{1}$ fixes $z_1 = \infty$ and maps $z_2\mapsto 0$. Henceforth, we may assume $z_2 = 0$, and further that $Z_2 = \minimatrix{0}{0}{0}{1}$.
Next, write $Z_3 = \minimatrix{e}{f}{\bar f}{h}$ with $eh = |f|^2$. Then $\langle Z_1, Z_3 \rangle = -\frac{h}{2}$ and $\langle Z_2, Z_3 \rangle = -\frac{e}{2}$. It follows that $e,h > 0$, so that $|f|^2 > 0$. As above, we have that $h \cdot Z_3 = \twovector{f}{h} \begin{pmatrix} \bar f & h \end{pmatrix}$, so that $z_3 = \frac{f}{h} \in \BB$ and $|z_3|^2 > 0$. Then $\bminimatrix{1}{0}{0}{z_3}$ defines an element of $\PGL^+(2,\BB)$ which fixes $z_1 = \infty$ and $z_2 = 0$ and maps $z_3 \mapsto 1$.
\end{proof}

\begin{Remark}
In the case $\mathcal B = \mathbb C$, $\mathbb X = \mathbb H^3$, the signs of the $Z_i$ can always be chosen to satisfy condition~\ref{negativity}. This is achieved by choosing the $Z_i$ to have positive trace.
\end{Remark}

\subsection{Shape parameters}

Let $Z_1, Z_2, Z_3, Z_4 \in \Herm(2,\BB)$ have rank one, and let $z_1, z_2, z_3, z_4$ denote the corresponding elements of $\mathbb P^1 \BB$. Assume that $Z_1, Z_2, Z_3$ determine an ideal triangle in $\XX$. By Proposition~\ref{prop:standard-position}, there is a unique $A \in \PGL^+(2,\BB)$ such that $A z_1 = \infty,  Az_2 = 0$, and $A z_3 = 1$. Then $$(z_1,z_2; z_3,z_4) := A z_4$$ is an invariant of the ordered ideal points $z_1,\ldots,z_4$, generalizing the cross ratio for $\CP^1$. 

\begin{Proposition}
$z_1,z_2,z_3,z_4$ define an ideal tetrahedron in $\mathbb X$ if and only if $z = (z_1,z_2;z_3,z_4)$ lies in $\mathcal B \subset \mathbb P^1 \mathcal B$ and satisfies:
\begin{equation}
|z|^2, |1-z|^2 > 0.
\label{spacelike}
\end{equation}
In this case $z$ is called the \emph{shape parameter} of the ideal tetrahedron (with ordered vertices).
\label{spacelike-prop}
\end{Proposition}
\begin{proof}
Assume the $z_i$ are in standard position, and choose representatives
\begin{align*}
Z_1 &= \btwovector{1}{0}\browvector{1}{0} = \bminimatrix{1}{0}{0}{0}&
Z_2 &= \btwovector{0}{1}\browvector{0}{1} = \bminimatrix{0}{0}{0}{1}\\
Z_3 &= \btwovector{1}{1}\browvector{1}{1} = \bminimatrix{1}{1}{1}{1}&
Z_4 &= \btwovector{a}{b}\browvector{\bar a}{\bar b} = \bminimatrix{|a|^2}{a\bar b}{\bar a b}{|b|^2}.
\end{align*}
We are free to change the signs of the $Z_i$ (or even multiply by a non-zero real number) with the aim that $Z_i$ satisfy Condition~(\ref{negativity}). First, note that $\langle Z_1, Z_2 \rangle, \langle Z_1, Z_3 \rangle, \langle Z_2, Z_3 \rangle = -\frac{1}{2} < 0$, so it will not be fruitful to change the signs of $Z_1,Z_2$ or $Z_3$. Next, 
\begin{align*}
\langle Z_1, Z_4 \rangle &= -\frac{|b|^2}{2}&
\langle Z_2, Z_4 \rangle &= -\frac{|a|^2}{2}&
\langle Z_3, Z_4 \rangle &= -\frac{|a-b|^2}{2}.
\end{align*}
Condition~(\ref{negativity}) is satisfied if and only if $|a|^2, |b|^2$ and $|a-b|^2 > 0$, which is true if and only if  $\btwovector{a}{b} \sim \btwovector{z}{1}$ with $|z|^2, |z-1|^2 > 0$.
\end{proof}

Using the language of Lorentzian geometry, we say that $z$ and $z-1$, as in the Proposition, are \emph{space-like}. In fact, all facets of an ideal tetrahedron are space-like and totally geodesic with respect to the metric induced by $\langle \cdot, \cdot \rangle$ on $\mathbb X$.

\smallskip

{\noindent \bf Orientation.} Given ideal vertices $z_1, z_2, z_3, z_4$, recall the map $i: \Delta^3 \to \XX$ from the standard three-dimensional ideal simplex $\Delta^3$ to the geometric ideal tetrahedron $T = T(z_1, z_2, z_3, z_4)$ determined by those vertices. The ideal tetrahedron is called \emph{positively oriented} if this map takes the standard orientation on $\Delta^3$ to the orientation of $\RP^3$ determined by the standard ordered basis $\bminimatrix{1}{0}{0}{1}$, $\bminimatrix{1}{0}{0}{-1}$, $\bminimatrix{0}{1}{1}{0}$, $\bminimatrix{0}{\kappa}{-\kappa}{0}$ coming form the coordinates~(\ref{coordinates-on-Herm}). This condition can be read off from the shape parameter. Orient $\partial \XX$ so that the direction pointing away from $\XX$ is positive. Then $\mathcal B \subset \partial \mathbb X$ inherits an orientation allowing us to make sense of the notion of positive imaginary part (use for example the ordered basis $\{ 1, \kappa\}$). The following is straightforward.

\begin{Proposition}
\label{prop:positively-oriented}
The ideal tetrahedron $T$ determined by $z_1, z_2, z_3, z_4$ is positively oriented if and only if its shape parameter $z$ has positive imaginary part.
\end{Proposition}

Recall that the ideal tetrahedron $T$ determined by $z_1, z_2, z_3, z_4$ is called a \emph{collapsed} $\mathbb H^2$ tetrahedron if $T$ is contained in a hyperbolic plane. Note that $T$ is collapsed if and only if its shape parameter $z$ is real.

The shape parameter $z = (z_1,z_2; z_3,z_4)$ is a natural geometric quantity associated to the edge $e = z_1z_2$ in the following sense, familiar from Thurston's notes in the hyperbolic case.  Change coordinates (using an element of $\PGL^+(2,\BB)$) so that $z_1 = \infty$, and $z_2 = 0$. Then the subgroup $G_e$ of $\PGL^+(2,\BB)$ that preserves $e$ is given by $$G_e = \left\{ A = \bminimatrix{\lambda}{0}{0}{1}: \lambda \in \BB, |\lambda|^2 > 0 \right\}.$$ The number $\lambda = \lambda(A)$ associated to $A \in G_e$ is called the \emph{exponential $\BB$-length} and generalizes the exponential complex translation length of a loxodromic element of $\PSL(2,\CC)$. Let $A \in G_e$ be the unique element so that $A z_3 = z_4$. Then the shape parameter is just the complex length of $A$: $z = \lambda(A)$. 

\begin{figure}[h]
{\centering

\def\svgwidth{1.5in}
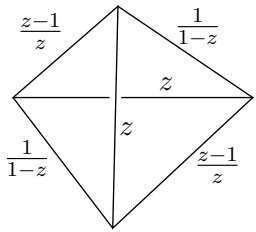

}
\caption{\label{ch5-shape-params} The shape parameters corresponding to the six edges of an ideal tetrahedron.}
\end{figure}

There are shape parameters associated to the other edges as well. We may calculate them as follows. Let $\pi$ be any \emph{even} permutation of $\{1,2,3,4\}$, which corresponds to an orientation preserving diffeomorphism of the simplex $\Delta^3$. Then $(z_{\pi(1)}, z_{\pi(2)}; z_{\pi(3)}, z_{\pi(4)})$ is the shape parameter associated to the edge $e' = z_{\pi(1)} z_{\pi(2)}$. This definition a priori depends on the orientation of the edge $e'$. However, we check that $$(z_2,z_1;z_4,z_3) = (z_1,z_2; z_3,z_4).$$
We need only demonstrate this in the case that the $z_i$ are in standard position. Consider the map $B = \minimatrix{0}{z}{1}{0}$, which exchanges $\infty$ and $0$ and also exchanges $1$ and $z$. Then 
\begin{equation*}
(0,\infty; z,1) = B\cdot 1 = z = (\infty,0;1,z).
\end{equation*}
Further, by considering $$B' = \bminimatrix{0}{1}{-1}{1} : \begin{array}{clc} \infty &\mapsto& 0 \\ 0 &\mapsto & 1 \\ 1 &\mapsto& \infty \end{array}$$ we obtain the relationship $$(1, \infty; 0, z) = B' \cdot z = \frac{1}{1-z}.$$ 
Figure~\ref{ch5-shape-params} summarizes the relationship between the shape parameters of the six edges of an ideal tetrahedron, familiar from the hyperbolic setting.

\subsection{Glueing tetrahedra together}
The faces of an ideal tetrahedron are hyperbolic ideal triangles. Given two tetrahedra $T,S$ and a face $\triangle v_1 v_2 v_4, \triangle w_1 w_2 w_3$ on each such that the orientations are opposite, there is a unique isometry $A \in \PGL^+(2,\mathcal B)$ mapping 
\begin{align*}
w_1 &\mapsto v_1,&
w_2 &\mapsto v_2,&
w_3 &\mapsto v_4
\end{align*}
that glues $S$ to $T$ along the given faces.
Suppose $T, S$ are in standard position so that $v_1 = w_1 = \infty$, $v_2 = w_2 = 0$, $v_3 = w_3 = 1$ and $v_4 = z$. Then the glueing map $A$ fixes $\infty$ and $0$ and acts as a (linear) similarity of $(\mathcal B, |\cdot|^2)$. This similarity is exactly multiplication by the shape parameter $z$ associated to the edge $v_1v_2$ of $T$. 
Composition of glueing maps for tetrahedra $T_1, T_2, \ldots, T_n$ in standard position about the common edge $0 \infty$ is described by the product of shape parameters $z_1 z_2 \ldots z_{n-1}$ (this describes the map that glues $T_n$ on to the other $n-1$ tetrahedra which have already been glued together).
Hence, in order for the geometric structure to extend over an interior edge $e$ of a union of tetrahedra, our shape parameters must satisfy:
\begin{eqnarray}
\prod_{T_i \text{ meets edge } e} z_i &=& 1 \label{edge}
\end{eqnarray}
where $z_i$ is the shape parameter associated to the edge of $T_i$ being identified to $e$. In fact we need that the development of the tetrahedra around the edge $e$ winds around the edge exactly once (in other words $\prod z_i$ is a rotation by $2\pi$ rather than $2\pi n$ for some $n \neq 1$). The terminology we will use for this condition is the following:

\begin{Definition}
\label{2pi}
Consider the $\mathbb X$ structure obtained by gluing together the ideal tetrahedra $T_1,\ldots,T_n$ in sequence around a single edge, denoted~$e$. The degree to which this structure is singular along~$e$ is measured by the holonomy of the development of the tetrahedra going once around $e$; here the holonomy is an isometry of $\XX$ that preserves the image of $e$ in $\XX$. We say that $e$ has \emph{total dihedral angle $2\pi$} if this edge holonomy has rotational part exactly~$2\pi$. 
\end{Definition}
\noindent Note that in some cases, $\mathbb X$ will not have a continuous group of isometries that rotate around a geodesic line. Nonetheless, rotations by multiples of $\pi$ around any geodesic are always defined in $\mathbb X$. See \cite{Danciger-13-1} for an explicit description of the local isometry group at a geodesic in each of the three geometries. 

\begin{figure}[h]
{\centering

\def\svgwidth{4.0in}
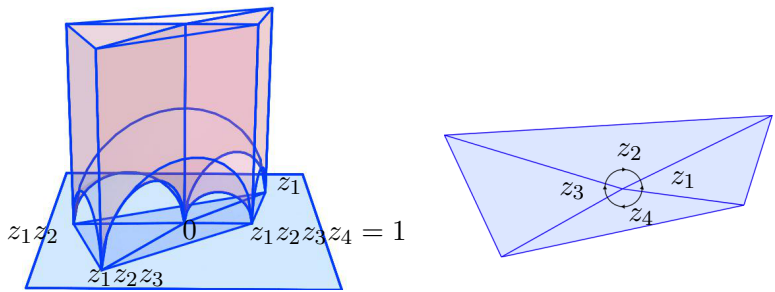 

}
\caption{The shape parameters going around an edge must have product one and total dihedral angle $2\pi$.}
\end{figure}

Let $M$ be a three-manifold with a fixed topological ideal triangulation $\mathcal{T} = \{\mathcal T_1,\ldots,\mathcal T_n\}$, that is $M$ is the union of tetrahedra $\mathcal T_i$ glued together along faces, with vertices removed. A triangulated $\mathbb X$ structure on $M$ is a realization of all the tetrahedra comprising $M$ as geometric tetrahedra so that the structure extends over all interior edges of the triangulation. This amounts to assigning each tetrahedron $\mathcal T_i$ a shape parameter $z_i$, such that for each interior edge $e$, the equation~(\ref{edge}) holds and $e$ has total dihedral angle $2\pi$. All of these equations together make up \emph{Thurston's equations} (also commonly called the \emph{edge consistency equations}). The solutions of these equations make up the \emph{deformation variety} of triangulated $\mathbb X$ structures on $M$.

%

\section{Triangulated geometric structures}
\label{triangulated-structures}
We apply the general construction from the previous section to build triangulated geometric structures for the cases $\mathbb X = \mathbb H^3, \mathbb H^2, \AdS^3, \HP^3$. Throughout, let $M$ be a three-manifold with a union of tori as boundary, and let $\mathcal{T} = \{\mathcal T_1,\ldots,\mathcal T_n\}$ be a fixed topological ideal triangulation of $M$.

\subsection{Triangulated $\HH^3$ structures: $\BB = \CC$}
Let $\kappa^2 = -1$, so that $\mathcal B = \mathbb C$ is the complex numbers. In this case, the inner product $\langle \cdot, \cdot \rangle$ on $\Herm(2,\mathbb C)$ is of type $(3,1)$ and $\mathbb X$ is the projective model for $\mathbb H^3$. Since $|z|^2 \geq 0$ holds for any $z$, with equality if and only if $z = 0$, Proposition~\ref{spacelike-prop} gives the well-known fact that any $z \in \mathbb C\setminus \{0,1\}$ is a valid shape parameter defining an ideal tetrahedron in $\mathbb H^3$.

Thus, hyperbolic structures on $(M, \mathcal{T})$ are obtained by solving Thurston's equations~(\ref{edge}) over $\mathbb C$ with all shape parameters $z_i$ having positive imaginary part.

\begin{Example}{(Figure eight knot complement)}
\label{fig8-example}
Let $M$ be the figure eight knot complement. Let $\mathcal T$ be the decomposition of $M$ into two ideal tetrahedra (four faces, two edges, and one ideal vertex) well-known from~\cite{Thurston}. 
\begin{figure}[h]
\center
\includegraphics[height = 1.2in] {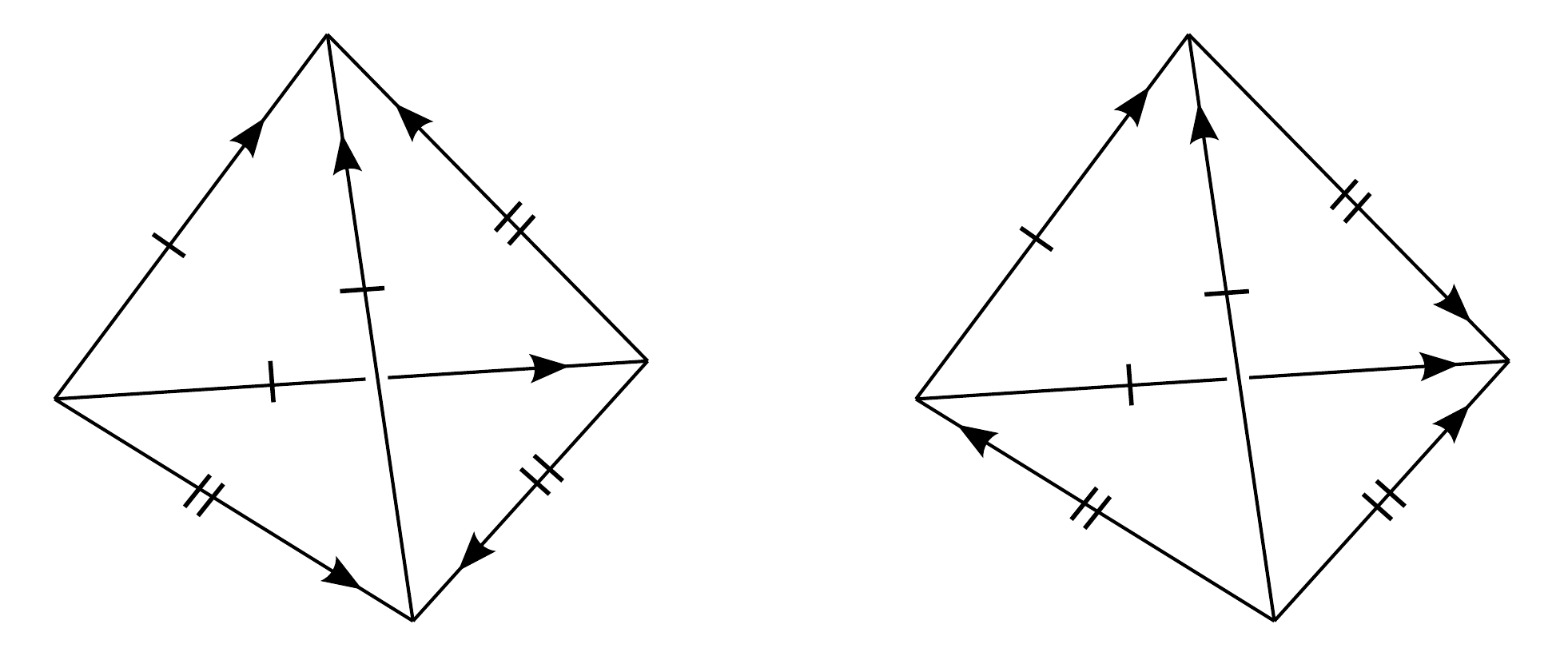}
\caption[The figure eight knot complement is the union of two ideal tetrahedra.]{The figure eight knot complement is the union of two ideal tetrahedra. In the diagram, identify two faces if the boundary edges and their orientations match. \label{fig8-2tets}}
\end{figure}
The edge consistency equations reduce to the following:
\begin{equation}z_1(1-z_1)z_2(1-z_2) = 1. \label{fig8-edge1}\end{equation}
The exponential complex length of the longitude $\ell$ and the meridian $m$, which can be read off from the triangulation of $\partial M$ (see Figure~\ref{boundary-fig8}), are given by
 \begin{align*}
H(\ell) &= z_1^2(1-z_1)^2 &
H(m) &= z_2(1-z_1).
\end{align*}
\begin{figure}[h]
{\centering

\def\svgwidth{4.0in}
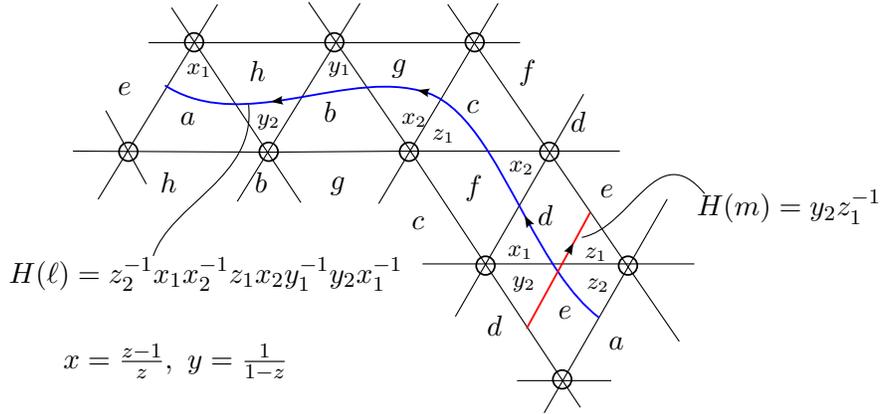

}
\caption[The boundary tessellation of the figure eight complement]{\label{boundary-fig8}The exponential complex lengths of $\ell$ and $m$ can be read off from a picture of the tessellation of $\partial M$. The triangles are labeled as in \cite[Ch.~4]{Thurston}}
\end{figure}
Let $\theta \in (0,2\pi)$ and enforce the additional condition that $$H(\ell) = z_1^2 (1-z_1)^{2} = e^{i\theta}.$$
A positively oriented solution is given by 
\begin{align*}
z_1 &= \frac{1 \pm \sqrt{1 - 4e^{i\theta /2}}}{2}& z_2 &= \frac{1 \pm \sqrt{1 - 4e^{-i\theta /2}}}{2}
\end{align*} where we may choose the root with positive imaginary part. 
The solution gives a hyperbolic structure whose completion $\bar M$ is topologically the manifold $M_{\ell}$ gotten by Dehn filling $M$ along $\ell$. The completed hyperbolic structure has a \emph{cone singularity} with cone angle $\theta$ (see e.g. \cite{Hodgson-05}). Note that $\bar M$ is a torus bundle over the circle with monodromy $\minimatrix{2}{1}{1}{1}$ and the singular locus is a curve running once around the circle direction.  

\end{Example}

\subsection{Transversely hyperbolic foliations: $\BB = \RR$}
\label{flattened-tets}
Consider the degenerate case $\mathcal B = \mathbb R$. Then $\Herm(2,\mathbb R)$ is the symmetric real matrices (which is $\mathbb R^3$ as a vector space) and $\langle \cdot, \cdot \rangle$ is of signature $(2,1)$. The resulting geometry is $\mathbb X = \mathbb H^2$. Proposition~\ref{spacelike-prop} gives that any $z \in \mathbb R \setminus \{ 0,1\}$ is a valid shape parameter defining an ideal tetrahedron in $\mathbb H^2$. Such tetrahedra are \emph{collapsed} (see the discussion at the beginning of Section~\ref{general-tets}). 

\begin{Proposition}
A solution to Thurston's equations~(\ref{edge}) over $\mathbb R$ defines a transversely hyperbolic foliation on $M$. Such a structure will be referred to as a \emph{triangulated transversely hyperbolic foliation} on $(M,\mathcal T)$. The deformation variety $\mathscr D_\RR$ of these structures is called the \emph{real deformation variety}.
\end{Proposition}
\begin{proof}
The real shape parameter $z_j \in \mathbb R$ assigned to the tetrahedron $\mathcal T_j$, determines a submersion from $\mathcal T_j$ onto an ideal quadrilateral in $\mathbb H^2$ sending faces of $\mathcal T_j$ to ideal triangles. Beginning with one base tetrahedron, these submersions can be developed to produce a globally defined local submersion $$D : \widetilde M \rightarrow \mathbb H^2$$ which is equivariant with respect to a representation $$\rho : \pi_1 M \rightarrow \PGL(2,\mathbb R).$$
That the shape parameters satisfy Thurston's equations guarantees that the map $D$ can be extended, still as a local submersion, over the interior edges of the triangulation. The map $D$ is a degenerate developing map defining a transversely hyperbolic foliation on $M$.
\end{proof}

Note that, in this case, the condition that the development of tetrahedra around an interior edge have total dihedral angle $2\pi$ is equivalent to requiring that exactly two of the $z_i$ at that edge be negative. An edge with negative shape parameter is thought of as having dihedral angle $\pi$, while an edge having positive shape parameter has dihedral angle zero.

\begin{Remark}
In the case of non-positively oriented solutions to Thurston's equations over $\CC$ (which do not directly determine $\HH^3$ structures), it is possible for the dihedral angle at an edge of some tetrahedron, defined via analytic continuation, to lie outside the range $(0,\pi)$. In the case that a path of such solutions converges to a real solution, each dihedral angle converges to $k\pi$ for some $k$, possibly with $k \neq 0,1$. We ignore these real solutions; there are no positively oriented solutions nearby.
\end{Remark}

\begin{Example}{(Figure eight knot complement)}
\label{fig8-foliations}
Let $M$ be the complement of the figure eight knot as defined in Example~(\ref{fig8-example}). To find transversely hyperbolic foliations on $M$, we solve the edge consistency equations 
\begin{equation}z_1(1-z_1)z_2(1-z_2) = 1. \label{fig8-edge2}\end{equation}
over $\mathbb R$. The variety of solutions to (\ref{fig8-edge2}) has four (topological) components:
\begin{align*}
 1. \ \ z_1 < 0 &\text{ and } z_2 <0 &  2. \ \ z_1 < 0 &\text{ and } z_2 >1\\
3. \ \ z_1 >1 &\text{ and } z_2 < 0 & 4. \ \ z_1 > 1 &\text{ and } z_2 > 1
\end{align*}

Cases 1 and 4 determine solutions with angular holonomy $4\pi$ around one edge and zero around the other edge. So these solutions are discarded. Cases 2 and 3 are symmetric under switching $z_1$ and $z_2$. So, the transversely hyperbolic structures on $(M,\mathcal T)$ are parametrized by $z_1 < 0$ (which determines $z_2 > 1$). It follows that the structures are also parametrized by $H(\ell) = z_1^2(1-z_1)^2$. This is a special case of Theorem~\ref{punctured-torus}.
\end{Example}

\subsection{Triangulated $\AdS^3$ structures and the pseudo-complex numbers}
\label{AdS-tetrahedra}

Let $\mathcal{B}$ be the real algebra generated by an element $\tau$, with $\tau^2 = +1$. 
As a vector space $\mathcal{B} = \mathbb R + \mathbb R \tau$ is two dimensional over $\mathbb R$. 
In this case, the form $\langle \cdot, \cdot \rangle$ on $\Herm(2,\mathcal B)$ is of signature $(2,2)$ and $\mathbb X = \AdS^3$ is the \emph{anti de Sitter} space. Before constructing triangulated $\AdS$ structures, we discuss some important properties of the algebra $\mathcal B = \mathbb R + \mathbb R \tau$.

\bigskip

{\noindent \bf The algebra $\BB = \mathbb R + \mathbb R \tau$ of pseudo-complex numbers.}
First, note that $\mathcal{B}$ is not a field as e.g. $(1+\tau)\cdot(1-\tau) = 0.$
The square-norm defined by the conjugation operation $|a+b\tau|^2 = (a + b \tau) \overline{(a + b\tau)} = a^2 - b^2 ,$
 comes from the $(1,1)$ Minkowski inner product on $\mathbb R^2$ (with basis  $\{1, \tau\}$). The space-like elements of $\mathcal{B}$ (i.e. square-norm $ > 0$), acting by multiplication on $\mathcal{B}$ form a group and can be thought of as the similarities of the Minkowski plane that fix the origin. Note that if $|a + b\tau|^2 = 0$ then $b=\pm a$ and multiplication by $a + b\tau$ collapses all of $\mathcal{B}$ onto the light-like line spanned by $a + b\tau$.

The elements $\frac{1+\tau}{2}$ and $\frac{1-\tau}{2}$ are two spanning idempotents which annihilate one another: $$\left(\frac{1\pm\tau}{2}\right)^2 = \frac{1\pm\tau}{2}, \ \text{ and } \ \left(\frac{1+\tau}{2}\right) \cdot  \left(\frac{1-\tau}{2}\right) = 0.$$ Thus $\mathcal{B} \cong \mathbb R \oplus \mathbb R$ as $\mathbb R$ algebras via the isomorphism 
\begin{equation}
a \left(\frac{1+\tau}{2}\right) + b \left(\frac{1-\tau}{2}\right) \longmapsto (a,b).
\label{isomorphism}
\end{equation}
We have a similar splitting for the algebra of $2 \times 2$ matrices $M_2 (\mathcal B)$: $$\left(\frac{1+\tau}{2} A + \frac{1-\tau}{2}B\right)\cdot \left(\frac{1+\tau}{2} C + \frac{1-\tau}{2}D\right) = \left(\frac{1+\tau}{2} AC + \frac{1-\tau}{2}BD\right)$$ and also $$\text{det}\left(\frac{1+\tau}{2} A + \frac{1-\tau}{2}B\right) = \frac{1+\tau}{2}\text{det}(A) + \frac{1-\tau}{2}\text{det}(B). $$
The orientation preserving isometries $\Isom^+ \AdS^3 = \PGL^+(2,\mathcal B)$ correspond to the subgroup of $\PGL(2,\mathbb R) \times \PGL(2,\mathbb R)$ such that the determinant has the same sign in both factors. The identity component of the isometry group (which also preserves time orientation) is given by $\PSL(2,\RR) \times \PSL(2,\RR)$.

\begin{Proposition}\label{prop:product-structure}
There is a natural isomorphism $\mathbb P^1 \mathcal B \cong \mathbb P^1 \mathbb R \times \mathbb P^1 \mathbb R$ that identifies the action of $\PGL(2, \mathcal B) $ on $\mathbb P^1 \mathcal B$ with that of $\PGL(2,\mathbb R) \times \PGL(2, \mathbb R)$ on $\RP^1 \times \RP^1$.
\end{Proposition}

\begin{proof}
The isomorphism $\RP^1 \times \RP^1 \rightarrow \mathbb P^1 \mathcal B$ is given by $$ \btwovector{a}{b}, \btwovector{c}{d} \longmapsto \frac{1+\tau}{2} \btwovector{a}{b} + \frac{1-\tau}{2} \btwovector{c}{d}.$$
\end{proof}

\noindent The space $\mathbb P^1 \mathcal B$, which is the \emph{Lorentz compactification} of $\mathcal B$, is covered by two copies of $\mathcal B$, the standard copy $\mathcal B \cong \left\{ \btwovector{x}{1} : \ x \in \mathcal B \right\}$ and another copy $\mathcal B \cong \left\{ \btwovector{1}{x} : \ x \in \mathcal B \right\}$. 
The square-norm $|\cdot|^2$ on $\mathcal B$ induces a flat conformal Lorentzian structure on each of these charts which agrees on the overlap, therefore defining a flat conformal Lorentzian structure on $\mathbb P^1 \mathcal B$. It is simple to check that this structure is preserved by $\PGL^+(2, \mathcal B)$. We refer to $\PGL^+(2, \mathcal B)$ as the \emph{Lorentz Mobius transformations}. With its conformal structure $\mathbb P^1 \mathcal B$ is the $(1+1)$-dimensional Einstein universe $\text{Ein}^{1,1}$ (see e.g. \cite{Barbot-08}). Note that in this $1+1$ dimensional setting, the conformal Lorentzian structure on  $\mathbb P^1 \mathcal B$ is entirely determined by the two null directions in each tangent space; these directions are exactly the coordinate directions in the product structure given by Proposition~\ref{prop:product-structure}.

\bigskip

{\noindent \bf Thurston's equations for $\AdS^3$.} 
We think of $\mathbb R + \mathbb R \tau$ as the Lorentzian plane equipped with the metric induced by $|\cdot |^2$. Proposition~\ref{spacelike-prop} immediately implies:
\begin{Proposition}
The following are equivalent:
\begin{enumerate}
\item The ideal vertices $z_1,z_2,z_3,z_4$ determine an ideal tetrahedron $T$.
\item The shape parameter $z = (z_1,z_2; z_3,z_4)$ of the edge $z_1 z_2$ is defined and satisfies $|z|^2 , |1-z|^2 > 0$.
\item The shape parameters $z, \frac{1}{1-z}, \frac{z-1}{z}$ of all edges of $T$ are defined and space-like.
\item Placing $z_1$ at $\infty$, the triangle $\triangle z_2 z_3 z_4$ has space-like edges in the Minkowski plane $\BB$.
\end{enumerate}
\end{Proposition}
\begin{figure}[h]
{\centering

\def\svgwidth{5.2in}
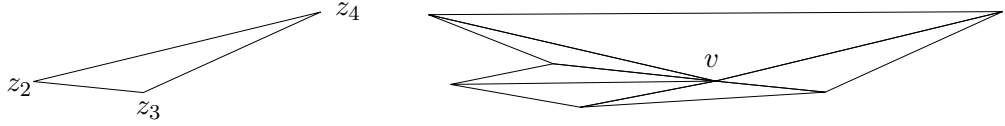

\vspace{0.1in}
}
\caption[The ``shadow" of an $\AdS$ ideal tetrahedron on the Lorentz plane.]{left: Placing one vertex $z_1$ at infinity, the other three vertices $z_2, z_3,z_4$ determine a \emph{spacelike} triangle in the Lorentzian plane: $|z_2 - z_3|^2, |z_3-z_4|^2,|z_4-z_2|^2 > 0$. right: if tetrahedra are glued together along an interior edge (connecting $\infty$ to $v$), the corresponding space-like triangles must fit together around the vertex $v$.}
\end{figure}

Similar to the case of degenerate tetrahedra, the total dihedral angle condition is discrete for $\AdS$ tetrahedra: 
\begin{Proposition}
The condition of Definition~\ref{2pi}, that the total dihedral angle around an interior edge be $2\pi$, is equivalent to the condition that exactly two of the $z_i$ at that edge have negative real part. 
\end{Proposition}

Using the isomorphism~(\ref{isomorphism}), a shape parameter $z \in \mathbb R + \mathbb R \tau$ can be described as a pair $(\lambda, \mu)$ of real numbers: $$z = \frac{1+\tau}{2} \lambda + \frac{1-\tau}{2} \mu.$$
Observe that \begin{align*}
|z|^2 &= \lambda \mu, &
|1-z|^2 &= (1-\lambda)(1-\mu).
\end{align*}
Hence $|z|^2, |1-z|^2 > 0 $ if and only if $\lambda$ and $\mu$ have the same sign and $1-\lambda$ and $1-\mu$ have the same sign. Thus, by Proposition~\ref{spacelike-prop},
$z$ is the shape parameter for an ideal tetrahedron in $\AdS$ if and only if $\lambda, \mu$ lie in the same component of $\mathbb R \setminus \{0,1\}$.

Further, as the imaginary part of $z$ is $\frac{\lambda - \mu}{2}$, a tetrahedron with shape parameter $z$ is positively oriented if and only $\lambda > \mu$ by Proposition~\ref{prop:positively-oriented}.

These observations combine to give:
\begin{Proposition}
The shape parameters $z_i =  \frac{1+\tau}{2} \lambda_i + \frac{1-\tau}{2} \mu_i$, for $i=1,\ldots,n$, define positively oriented ideal tetrahedra that glue together compatibly around an edge in $\AdS^3$ if and only if: 
\begin{itemize}\item $\prod_{i=1}^n \lambda_i = 1 \text{\ \ \  and \ \ \ } \prod_{i=1}^n \mu_i = 1$.
\item $\lambda_i,\mu_i$ lie in the same component of $\mathbb R \setminus \{0,1\}$ for each $i = 1,\ldots,n$, 
\item $\lambda_i > \mu_i$ for each $i=1, \ldots, n$.
\item For exactly two $i \in \{1,\ldots,n\}$, we have $\lambda_i,\mu_i < 0$.
\end{itemize}
\label{AdS-RcrossR}
\end{Proposition}

Thus a triangulated $\AdS$ structure on $(M, \mathcal T)$ is determined by two triangulated transversely hyperbolic foliations on $(M,\mathcal T)$ whose shape parameters $(\lambda_i)$ and $(\mu_i)$ obey the conditions set out in the above Propositions. This gives a concrete method for regenerating $\AdS$ structures from transversely hyperbolic foliations:

\begin{Corollary}
Let $(\lambda_i)$ be shape parameters defining a transversely hyperbolic foliation on $(M,\mathcal T)$. Suppose this structure can be deformed to a new one with shape parameters $\lambda_i' > \lambda_i$. Then $z_i = \frac{1+\tau}{2} \lambda_i + \frac{1-\tau}{2} \lambda_i'$ defines an $\AdS$ structure on $(M, \mathcal T)$.
\end{Corollary}
\noindent This suggests the following more general question, which we do not address here:
\begin{Question}
When and how do two transversely hyperbolic foliations on $M$ determine an $\AdS$ structure in the absence of an ideal triangulation?
\end{Question}

\bigskip

{\noindent \bf Tachyons.}
Consider a triangulated manifold $(M,\mathcal T)$. Let us assume that there is only one ideal vertex $v$ in $\mathcal T$ (after identification). Then $\partial M$, which is naturally identified with $L(v)$, has only one component. Assume that $\partial M$ is a torus and that $M$ has a fixed $\AdS$ structure determined by a positively oriented solution to Thurston's equations over $\BB = \Rtau$. Let $\mathcal N (v)$ be a neighborhood in $M$ of the ideal vertex $v$. Similar to the hyperbolic case, the $\AdS$ structure on $M$ induces a structure on $L(v)$ modeled on the similarities of the Minkowski plane $\RR^{1,1}$ which we identify with $\BB$. The similarities of $\BB$ that fix the origin are exactly the space-like elements $\BB^+$ (i.e. the elements with positive square-norm) acting by multiplication. Just as in the hyperbolic case, the geodesic completion of $\mathcal N(v)$ can be understood in terms of this similarity structure. Let $D_{\partial}: \widetilde{L(v)} \rightarrow \BB$ be the developing map. Assuming that $M$ is not complete near $v$, the holonomy $H$ of the similarity structure on $L(v)$ fixes a point, which we may assume to be the origin. Then $H: \pi_1 L(v) \rightarrow \BB^+$ is the \emph{exponential $\BB$-length} function restricted to $\pi_1 \partial M$. The image of $D_\partial$ does not contain the origin, so $D_\partial$ determines a lift $\widetilde H$ of $H$ to the similarities $\widetilde \BB^+$ of $\widetilde{\BB \setminus 0}$: 
\begin{align*}
\widetilde H(\gamma) &= \log |H(\gamma)| + \tau \phi(\gamma) + i\tilde{ \mathcal R}(\gamma) \ \ \in \ \ \RR + \tau \RR + i\pi \ZZ
\end{align*}
where $\log |H(\gamma)|$ is the translation length of $H(\gamma)$, $\phi(\gamma)$ is the hyperbolic angle of the boost part of $H(\gamma)$ and $\tilde{ \mathcal R}(\gamma)$ is the \emph{total rotational part} of the holonomy of $\gamma$  (see Definition~7 of \cite{Danciger-13-1}), which is an integer multiple of $\pi$ measuring the number of half rotations around $0$ swept out by developing along $\gamma$. Assume that there is some element of $\pi_1 \partial M$ with non-zero discrete rotational part. Then $D_\partial$ is a covering map onto $\BB \setminus 0$ (this follows from the more general theory of affine structures on the torus \cite{Nagano-74}). In the half-space model for $\AdS$ (see Appendix~A of \cite{Danciger-11}), the developing map $D_{\partial}$ is the ``shadow" of the developing map $D : \widetilde{\mathcal N (v)} \rightarrow \AdS^3$. So the image of $D$ is $I \setminus \Cln$, where $I$ is a neighborhood of the geodesic $\Cln$ with endpoints $0, \infty$ (note that $I$ is not a cone as it is in the hyperbolic case). The completion of $\widetilde{\mathcal N(v)}$ is then given by adjoining a copy of $\Cln$. So the completion of $\mathcal N(v)$ is 
\begin{align*}
\overline{\mathcal N(v)} &= \left(\widetilde{I \setminus \Cln} \cup \Cln \right)/ \widetilde H(\pi_1 \partial M) = \mathcal N(v) \cup (\Cln / H(\pi_1 \partial M)).
\end{align*}

In particular, if the moduli $|H(\pi_1 \partial M)|$ form a discrete subgroup of the multiplicative group $\RR^+$, then $\Cln / H(\pi_1 \partial M)$ is a circle and $\overline{\mathcal N(v)}$ is a manifold. This is the case if and only if there exists a generator $\alpha$ of $\pi_1 \partial M$ such that $\widetilde H(\alpha)$ is a rotation by $k\pi \neq 0$ plus a boost. In this case, the completion $\overline M$ (which is given near the boundary by $\overline{\mathcal N(v)}$) is topologically the manifold $M_\alpha$ obtained by Dehn filling $M$ along the curve $\alpha$. If the discrete rotational part $\tilde{\mathcal R}(\alpha) = 2\pi$, then $\overline{M}$ has a \emph{tachyon} singularity (\cite{Barbot-09} or see Section~2.5 of \cite{Danciger-13-1}) with mass equal to the hyperbolic angle $\phi(\alpha)$ of the boost part of $H(\alpha)$. 

\begin{Remark}
If the moduli $|H(\pi_1 \partial M)|$ are dense in $\RR^+$, the geodesic completion of $M$ near $v$ has a topological singularity that resembles a Dehn surgery type singularity in hyperbolic geometry. This more general singularity has not yet been studied to the knowledge of the author. \end{Remark}

\begin{Example}{(figure eight knot complement)}
\label{fig8-AdS}
Let $M$ be the complement of the figure eight knot from Examples \ref{fig8-example} and \ref{fig8-foliations}.
We use Proposition~\ref{AdS-RcrossR} and the analysis in Example~\ref{fig8-foliations} to build $\AdS$ structures on $(M,\mathcal T)$. Consider the connected component of real solutions to the edge consistency Equation~(\ref{fig8-edge2}) with $z_1 < 0$ and $z_2 > 1$. Taking the differential of $\log$ of Equation~(\ref{fig8-edge2}), we obtain
\begin{equation}
\frac{2z_1}{z_1(1-z_1)}dz_1 + \frac{2z_2}{z_2(1-z_2)}dz_2 = 0
\end{equation}
which implies that $\frac{dz_2}{dz_1} > 0$ at any point of (this connected component of) the variety. Thus, any two distinct solutions $(\lambda_1,\lambda_2)$ and $(\mu_1,\mu_2)$ satisfy (up to switching the $\lambda$'s with the $\mu$'s) $$\lambda_1 > \mu_1 \text{ \ \ \ and \ \ \ } \lambda_2 > \mu_2$$ and give a positively oriented solution 
\begin{align*}
z_1 &= \frac{1+\tau}{2} \lambda_1 + \frac{1-\tau}{2} \mu_1 &
z_2 &= \frac{1+\tau}{2} \lambda_2 + \frac{1-\tau}{2} \mu_2
\end{align*}
to the edge consistency equations over $\mathbb R + \mathbb R\tau$ determining $\AdS$ structures on $M$. It is straight forward to show that the discrete rotational part of the holonomy of $\ell$ is $\tilde{\mathcal R}(\ell) = +2\pi$.

Now impose the additional condition \begin{eqnarray*}H(\ell) = z_1^2(1-z_1)^2 &=& e^{\tau \phi} := \cosh \phi + \tau \sinh \phi, \end{eqnarray*}
which is equivalent to \begin{eqnarray*} \lambda_1^2(1-\lambda_1)^2 &=& \cosh(\phi)+\sinh(\phi) \ = \ e^{\phi}\\ \mu_1^2(1-\mu_1)^2 &=& \cosh(\phi) - \sinh(\phi) \ = \ e^{-\phi}.\end{eqnarray*}
 
The geodesic completion of the $\AdS$ structure determined by $(z_1,z_2)$ is an $\AdS$ structure on the Dehn filled manifold $M_\ell$ with a tachyon of mass $\phi$. Note that as $\mu_1 < \lambda_1 < 0$, the tachyon mass is negative.
\end{Example}

\subsection{Triangulated $\HP$ structures}
\label{HP-tets}
 Next let $\mathcal B = \mathbb R + \mathbb R \sigma$ where $\sigma^2 = 0$.  Then the form $\langle \cdot, \cdot \rangle$ on $\Herm(2,\mathcal B)$ is degenerate (with the eigenvalue signs $+,+,-,0$). In this case, $\mathbb X = \HP^3$ and $\PGL^+(2,\mathbb R + \mathbb R \sigma) \cong G_{\HP}$ give the projective model for \emph{half-pipe geometry}.

We equip $\mathbb R + \mathbb R \sigma$ with the degenerate metric induced by $|\cdot |^2$. Proposition~\ref{spacelike-prop} immediately implies:
\begin{Proposition}
The following are equivalent:
\begin{enumerate}
\item The ideal vertices $z_1,z_2,z_3,z_4 \in \mathbb P^1 \Rsigma$ define an ideal tetrahedron $T$.
\item The shape parameter $z = (z_1,z_2; z_3,z_4)$ of the edge $z_1 z_2$ is defined and satisfies $\text{Re} \ z \neq 0,1$.
\item The shape parameters $z, \frac{1}{1-z}, \frac{z-1}{z}$ of all edges of $T$ are defined and have real parts not equal to $0,1$.
\item Placing $z_1$ at $\infty$, $\triangle z_2 z_3 z_4$ is a triangle in the $\mathbb R + \mathbb R \sigma$ plane that has non-degenerate edges.
\end{enumerate}
\end{Proposition}
The real part $a$ of $z = a + b\sigma$ describes a collapsed tetrahedron in $\mathbb H^2$, while the imaginary part $b\sigma$ describes an infinitesimal ``thickness". If $b > 0$, then the tetrahedron is positively oriented; In this case $z$ is thought of as being tangent to a path of complex (resp. $\mathbb R + \mathbb R \tau$) shape parameters describing a collapsing family of positively oriented hyperbolic (resp. $\AdS$) tetrahedra.

\begin{Proposition}
The shape parameters $z_i =  a_i + b_i \sigma$, for $i=1,\ldots,n$, define ideal tetrahedra that glue together compatibly around an edge in $\HP^3$ if and only if: 
\begin{itemize} \item $(a_1,\ldots,a_n) \in \mathbb R^n$ satisfy the equation $\prod_{i=1}^n a_i = 1$,
\item $(b_1,\ldots,b_n) \in T_a \mathbb R^n$ satisfy the differential of that equation \begin{eqnarray*} d(\prod_{i=1}^n z_i)\Big |_{z_i = a_i} (b_1,\ldots,b_n) = 0,\end{eqnarray*}
\item
and exactly two of the $a_i$ are negative. 
\end{itemize}
Thus the real part of a solution to Thurston's equations over $\mathbb R + \mathbb R \sigma$ defines a triangulated transversely hyperbolic foliation, and the imaginary ($\sigma$) part defines an infinitesimal deformation of this structure.  

\end{Proposition}

%
%
%

\section{Regeneration of $\mathbb H^3$ and $\AdS^3$ structures}
\label{s:tet-regen}
In the context of ideal triangulations, the problem of regenerating hyperbolic and $\AdS$ structures from transversely hyperbolic foliations becomes more straight-forward, especially in the presence of smoothness assumptions.

\begin{Proposition}
\label{tet-regen}
Let $(z_j) \in \RR^N$ be a solution to Thurston's equations determining a transversely hyperbolic foliation $\mathcal F$. Suppose the real deformation variety $\mathscr D_\RR$ is smooth at $(z_j)$, and suppose that $(v_j) \in \mathbb R^N$ is a tangent vector to $\mathscr D_\RR$ such that $v_j > 0$ for all~$j$. Then there are hyperbolic structures $\mathscr H_t$ and $\AdS$ structures $\mathscr A_t$ on $(M, \mathcal T)$, defined for $t > 0$, such that $\mathscr H_t$ and $\mathscr A_t$ collapse to $\mathcal F$ as $t \to 0$.
\end{Proposition}

\begin{proof}
The imaginary tangent vector $(iv_j)$ can be integrated to give a path of complex solutions to the edge consistency equations. Similarly, the imaginary tangent vector $(\tau v_j)$ can be integrated to give a path of $\mathbb R + \mathbb R \tau$ solutions. This is easy to see explicitly, for consider a smooth path $(\lambda_j(t))$ in $\mathscr D_\RR$, with $\lambda_j(0) = z_j$, and $\lambda_j'(0) = v_j$ at $t = 0$. Then, the path $$z_j(t) = \frac{1+\tau}{2} a_j(t) + \frac{1-\tau}{2} a_j(-t)$$ has $z_j'(0) = \tau v_j$.

 In both the hyperbolic and $\AdS$ cases the solutions have positive imaginary part, so they determine geometric structures.
\end{proof}

In light of this proposition, we ask the following question:
\begin{Question}
Given a triangulated three-manifold $(M,\mathcal T)$, which points of the real deformation variety $\mathscr D_\RR$ are smooth with positive tangent vectors? 
\end{Question}
\noindent Theorem~\ref{punctured-torus}, proved in Section~\ref{torus-bundles} gives a partial answer to this question in the case that $M$ is a punctured torus bundle. Theorem~\ref{thm:regen} follows from Theorem~\ref{punctured-torus} by Proposition~\ref{tet-regen}.

%
%
%

\section{Geometric transitions via triangulations}
\label{tet-transitions}

In this section, we show that when the conditions of Proposition~\ref{tet-regen} are satisfied, the regenerated hyperbolic and $\AdS$ structures may be organized into a path of real projective structures. The following proposition, in combination with Theorem~\ref{punctured-torus}, proves Theorem~\ref{triangulated-transition}.

\begin{Proposition}
\label{prop:tet-transition}
Suppose the real deformation variety $\mathscr D_\RR$ is smooth at a point $(z_j) \in \mathbb R^N$ and suppose $(v_j) \in \mathbb R^N$ is a positive tangent vector. Let $\mathscr H_t$ and $\mathscr A_t$ be hyperbolic $\AdS$ structures as in Proposition~\ref{tet-regen}. Then there is a path of projective structures $\mathscr P_t$ such that 
\begin{itemize}
\item $\mathscr P_t$ is conjugate to $\mathscr H_t$ for $t > 0$
\item $\mathscr P_t$ is conjugate to $\mathscr A_{|t|}$ for $t< 0$.
\item $\mathscr P_0$ is a half-pipe structure.
\end{itemize}
For each $t$, the triangulation $\mathcal T$ is realized by positively oriented tetrahedra in $\mathscr P_t$.
\end{Proposition}

The path of projective structures in the proposition is a \emph{geometric transition} from hyperbolic to $\AdS$ structures, as defined in \cite{Danciger-13-1}. In order to prove the proposition, we take the solutions to Thurston's equations given by Proposition~\ref{tet-regen} and realize them as a continuous path $(w_j(t))$ of solutions over shape parameter algebras $\BB_t$ which vary from $\BB_{+1} = \CC$ to $\BB_{-1} = \Rtau$.

\begin{proof}[Proof of Proposition~\ref{prop:tet-transition}]
Let $\BB_t = \RR + \RR \kappa_t$, where $\kappa_t^2 = -t|t|$. For $t > 0$, the map $\mathfrak a_t : \CC \to \BB_t$ defined by $i \mapsto \kappa_t/|t|$ is an isomorphism of algebras. For $t < 0$, the map $\mathfrak a_t : \Rtau \to \BB_t$ defined by $\tau \mapsto \kappa_t/|t|$ is an isomorphism.

Let $(\zeta_j(t)) \in \CC^N$  and $(\xi_j(t)) \in (\Rtau)^N$ be solutions to Thurston's equations with $\zeta_j(0) = \xi_j(0) = z_j$, $\zeta_j'(0) = i v_j$, $\xi_j'(0) = \tau v_j$, as guaranteed by Proposition~\ref{tet-regen}. Then, define $(w_j) \in \BB_t^N$ by:
\begin{itemize}
\item $w_j(t) = x_j(t) + \kappa_t y_j(t) := \mathfrak a_t (\zeta_j(t))$ if $t > 0$,
\item $w_j(t) = x_j(t) + \kappa_t y_j(t) :=  \mathfrak a_t (\xi_j(t))$ if $t < 0$.
\item $w_j(0) = x_j(0) + \kappa_0 y_j(0) := z_j + \kappa_0 v_j$.
\end{itemize}
The real part $x_j$ and the imaginary part $y_j$ of $w_j$ are continuous\marginnote{more regularity?} functions of $t$. Further, $y_j(t) > 0$ for all $t$ (in an open neighborhood of $t = 0$). Hence each solution $(w_j(t))$ determines a projective structure $\mathscr P_t$ built from positively oriented ideal tetrahedra in the model space $\XX_t$. As the functions $x_j, y_j$ vary continuously, as do the models $\XX_t$, we have that $\mathscr P_t$ varies continuously. Of course, for $t > 0$, the hyperbolic structure $\mathscr H_t$ is conjugate to the projective structure $\mathscr P_t$ via the obvious transformation, induced by $\mathfrak a_t$, conjugating $\HH^3 = \XX_{+1}$ to $\XX_t$ (see Section 4.5 of \cite{Danciger-13-1}). Similarly, for $t< 0$, $\mathscr A_{|t|}$ is conjugate to $\mathscr P_t$. 

Lastly, as $\kappa_0^2 = 0$, $\kappa_0$ is just another name for the element $\sigma$ in the description of half-pipe geometry (Section~\ref{models}). So $\mathscr P_0$ is a half-pipe structure.
\end{proof}

The shape parameters $w_j(t)$ from the proof of Proposition~\ref{prop:tet-transition} lie in different algebras $\BB_t$, making it slightly annoying to discuss continuity of the $w_j(t)$. One convenient solution to this issue is the following. Consider the generalized Clifford algebra $\mathcal C = \langle 1, i, \tau:\ i^2 = -1,\ \tau^2 = +1, \ i \tau = - \tau i \rangle$. We note that the algebras $\BB_t$, from the proof of Proposition~\ref{prop:tet-transition}, may be embedded as a continuous (even differentiable) path of sub-algebras inside $\mathcal C$ via the identification $\kappa_t = \left((1+t|t|)i + (1 - t |t|)\tau\right)/2$.
This allows us to think of $w_j(t)$ as a continuous path inside $\mathcal C$.

\begin{Example}
\label{fig8-ex}
Let $M$ be the figure eight knot complement, discussed in Examples~\ref{fig8-example}, \ref{fig8-foliations}, and \ref{fig8-AdS}. Let $\mathcal T$ be the decomposition of $M$ into two ideal tetrahedra (four faces, two edges, and one ideal vertex) well-known from~\cite{Thurston} (see Figure~\ref{fig8-2tets}). 
The edge consistency equations reduce to the following:
\begin{equation}z_1(1-z_1)z_2(1-z_2) = 1. \label{fig8-edge}\end{equation}
\noindent In Example~\ref{fig8-foliations}, we showed that the variety of real solutions to $(\ref{fig8-edge})$ (with total dihedral angle $2\pi$ around each edge) is a smooth one-dimensional variety with positive tangent vectors.
Thus, any transversely hyperbolic foliation on $(M,\mathcal T)$ regenerates to robust hyperbolic and $\AdS$ structures by Proposition~\ref{tet-regen}. As $M$ is a punctured torus bundle, this is a special case of Theorem~\ref{punctured-torus}.

Next, we consider hyperbolic cone structures on $M$, with singular meridian being the longitude $\ell$ of the knot (this is also the curve around the puncture in a torus fiber). Recall from Example~\ref{fig8-example} that such a structure, with cone angle $\theta < 2\pi$, is constructed by solving the equations
\begin{align}
\label{cone-equation}
H(\ell) = z_1^2(1-z_1)^2  &= e^{i\theta} = e^{-i(2\pi-\theta)}
\end{align}
over $\mathbb C$.
Recall from Example~\ref{fig8-AdS} that $\AdS$ tachyon structures with mass $\phi < 0$ are constructed by solving the equations
\begin{align}
 \label{tachyon-equation} 
H(\ell) = z_1^2(1-z_1)^2 &= e^{\tau\phi} = e^{-\tau(-\phi)} 
\end{align}
over $\mathbb R + \mathbb R \tau$. In order to construct a transition between these two types of structures, we consider a generalized version of these equations defined over the transitioning family $\mathcal B_t = \RR + \RR \kappa_t$ from the proof of Proposition~\ref{prop:tet-transition}. The idea is to replace $i$ in (\ref{cone-equation}) (resp. $\tau$ in (\ref{tachyon-equation})) by the algebraically equivalent elements $\kappa_t/|t|$. 
The generalized version of (\ref{cone-equation}) and (\ref{tachyon-equation}) that we wish to solve is 
\begin{align}
H(\ell) = z_1^2(1-z_1)^2 &= -e^{-\kappa_t}.
\label{generalized-equation}
\end{align}
\begin{figure}[h]
  \begin{center}
  \vspace{0.1in}
    \includegraphics[width=2.5in]{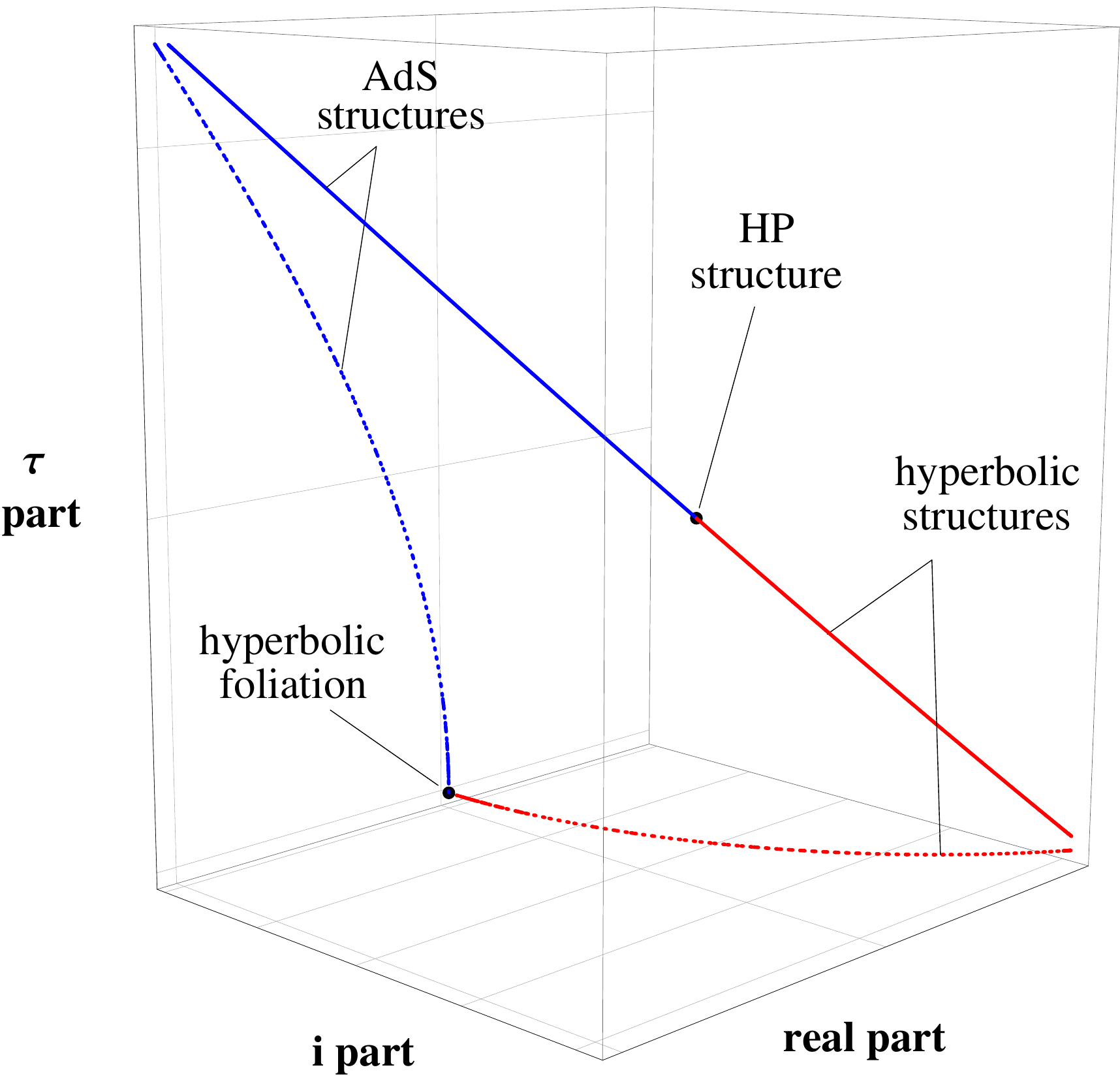}
  \end{center}
  \caption{Using the embeddings $\BB_t \hookrightarrow \mathcal C$ described above, we plot the $\mathcal{C}$-length of the singular curve as hyperbolic cone structures (red) transition to $\AdS$ tachyon structures (blue). After rescaling (solid lines), the transition is realized as a $\mathscr C^1$ path passing through a half-pipe structure.}
  \label{C-length}
\end{figure}
Note that the right hand side (which can be defined in terms of Taylor series) is a smooth function of $t$. In fact, solving (\ref{generalized-equation}) over the varying algebra $\mathcal B_t$ for small $t$, gives a differentiable path $(z_1(t),z_2(t))$ of shape parameters for transitioning structures. For $t > 0$, $(z_1,z_2)$ determines a hyperbolic cone structure with cone angle $\theta = 2\pi - |t|$. For $t < 0$, $(z_1,z_2)$ determines a $\AdS$ tachyon structure with hyperbolic angle $\phi = -|t|$.  At $t=0$, we relabel $\kappa_0 = \sigma$ for cosmetic purposes. Shape parameters defining the transitional half-pipe structure are given by: $$z_1(0) = \tfrac{1-\sqrt{5}}{2} + \tfrac{1}{2\sqrt{5}}\sigma, \ \  z_2(0) = \tfrac{1+\sqrt{5}}{2} + \tfrac{1}{2\sqrt{5}}\sigma.$$
The exponential $\Rsigma$-length of the curve $\ell$ around the singular locus is 
\begin{align*}
H(\ell) = z_1^2(1-z_1)^2 &= \left(\tfrac{1-\sqrt{5}}{2} + \tfrac{1}{2\sqrt{5}}\sigma\right)^2\left(\tfrac{1+\sqrt{5}}{2} - \tfrac{1}{2\sqrt{5}}\sigma\right)^2\\ &= (-1 + \tfrac{1}{2}\sigma)^2\\ &= 1 -\sigma \ = e^{-1\cdot\sigma}.
\end{align*}
The solution $(z_1(0),z_2(0))$ defines an $\HP$ structure whose completion has a cone-like singularity, called an infinitesimal cone singularity (see Section~4 of \cite{Danciger-13-1}). The infinitesimal cone angle is $-1$.
\end{Example}

%
%

\section{Punctured Torus Bundles}
\label{torus-bundles}

In this section we study the real deformation variety for $M$ a hyperbolic punctured torus bundle and prove Theorem~\ref{punctured-torus} from the introduction. 

\subsection{The monodromy triangulation}

We begin by describing the \emph{monodromy triangulation} (sometimes referred to as the Floyd-Hatcher triangulation) and Gueritaud's convenient description of Thurston's equations for this triangulation. 
See \cite{Gueritaud-06} for an elegant and self-contained introduction to this material. 

We think of the punctured torus $T^2$ as the quotient of $\mathbb R^2 \setminus \mathbb Z^2$ by the lattice of integer translations $\mathbb Z^2$. Any element of $\SL(2,\mathbb Z)$ acts on $T$ since it normalizes the lattice $\mathbb Z^2$. An element $\phi \in \SL(2,\mathbb Z)$ with distinct real eigenvalues $\lambda_+, \lambda_-$ is called \emph{Anosov}. We focus on the case that $\phi$ has positive eigenvalues. If $\phi$ has negative eigenvalues, then the following construction can be performed using $-\phi$ in place of $\phi$ with some small modifications; the resulting edge consistency equations will be the same. 
The following is a well-known fact (see e.g. \cite[Prop. 2.1]{Gueritaud-06} for a proof)
\begin{Proposition}
An Anosov $\phi \in \SL(2,\mathbb Z)$ with positive eigenvalues can be conjugated to have the following form:
\begin{equation*}
A \phi A^{-1} =: W  = R^{m_1}L^{n_1}R^{m_2}L^{n_2}\cdots R^{m_k}L^{m_k}
\end{equation*}
where $m_1,n_1 \ldots, m_k,n_k$ are positive integers and $R,L$ are the standard transvection matrices 
\begin{equation*}
R = \minimatrix{1}{1}{0}{1} \ \ \ \text{and } \ \ \ L = \minimatrix{1}{0}{1}{1}.
\end{equation*}
This form is unique up to cyclic permutation of the factors.
\end{Proposition}

This fact gives a canonical triangulation of the mapping torus $M = T \times I/ (x,0)\sim(\phi x, 1)$ as follows. Since $\phi$ and $W = A\phi A^{-1}$ produce homeomorphic mapping tori, we will henceforth assume $\phi = W$ has the form described in the proposition. We think of $W$ as a word of length $N = m_1+n_1+\ldots+m_k+n_k$ in the letters $L$ and $R$. Now, we begin with the standard ideal triangulation $\tau_0$ of $T$ having edges $(1,0), (0,1), (-1,1)$ (see figure below). Apply the first (left-most) letter of the word, which is $R$, to $\tau_0$ to get a new ideal triangulation $\tau_1 = R\tau_0$. These triangulations differ by a \emph{diagonal exchange}. Realize this diagonal exchange as an ideal tetrahedron as follows. Let $\mathcal T_1$ be an affine ideal tetrahedron in $T^2 \times \mathbb R$ with two bottom faces that project to the ideal triangles of $\tau_0$ in $T^2$ and two top faces that project to the ideal triangles of $\tau_1$ in $T^2$. 

\begin{figure}[!h]
{\centering

\def\svgwidth{2.5in}
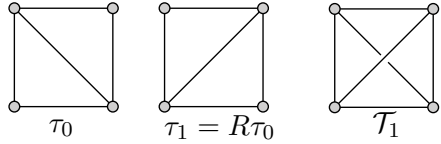

}
\caption{A diagonal exchange determines an ideal tetrahedron.}
\end{figure}
Next, we apply the first (left-most) two letters of $W$ to $\tau_0$ in order to get another ideal triangulation $\tau_2$. We note that $\tau_1$ and $\tau_2$ differ by a diagonal exchange and we let $\mathcal T_2$ be the ideal tetrahedron with bottom faces $\tau_1$ and top faces $\tau_2$. The bottom faces of $\mathcal T_2$ are glued to the top faces of $\mathcal T_1$. 
We proceed in this way to produce a sequence of $N+1$ ideal triangulations $\tau_0, \ldots, \tau_N$ with $\tau_k = W_k \tau_0$, where $W_k$ are the first (left-most) $k$ letters of $W$. It is easy to see that $\tau_{k}$ and $\tau_{k+1}$ differ by a diagonal exchange: for example if $W_{k+1} = W_k R$, then $W_{k+1}\tau_0$ and $W_k \tau_0$ differ by a diagonal exchange because $R \tau_0$ and $\tau_0$ differ by a diagonal exchange.  For consecutive $\tau_k, \tau_{k+1}$ define a tetrahedron $\mathcal T_{k+1}$ which has $\tau_k$ as its bottom faces and $\tau_{k+1}$ as its top faces. $\mathcal T_{k+1}$ is glued to $\mathcal T_{k}$ along $\mathcal \tau_k$. Note that the top ideal triangulation $\tau_N$ of the top tetrahedron $\mathcal T_{N}$ is given exactly by $\tau_N = \phi \tau_0$. So we glue $\mathcal T_{N}$ along its top faces $\tau_N$ to $\mathcal T_1$ along its bottom faces $\tau_0$ using the Anosov map $\phi$. The resulting manifold is readily seen to be $M$, the mapping torus of $\phi$. This decomposition into ideal tetrahedra is called the \emph{monodromy triangulation} or the \emph{monodromy tetrahedralization}.  

We note that the ideal triangulation $\tau_k$ of $T^2$ is naturally realized as a pleated surface inside $M$, at which the tetrahedra $\mathcal T_k$ and $\mathcal T_{k+1}$ are glued together. Further, we may label each $\tau_k$ with the $k^{th}$ letter of $W$. Hence, each tetrahedron $\mathcal T_{k+1}$ can be labeled with two letters, the letter corresponding to its bottom pleated surface $\tau_k$ followed by the letter corresponding to its top pleated surface $\tau_{k+1}$. If $\mathcal T_{k}$ is labeled $RL$ or $LR$ it is called a \emph{hinge} tetrahedron. Consecutive $LL$ tetrahedra make up an $LL$-\emph{fan}, while consecutive $RR$ tetrahedra make up an $RR$-\emph{fan}.

In order to build geometric structures using the monodromy triangulation, we assign shape parameters to the edges of the tetrahedra as follows: For tetrahedron $\mathcal T_i$, we assign the shape parameter $z_i$ to the (opposite) edges corresponding to the diagonal exchange taking $\tau_i$ to $\tau_{i+1}$. The shape parameters $x_i = \frac{z_i - 1}{z_i}$ and $y_i = \frac{1}{1-z_i}$ are assigned to the other edges according to the orientation of the tetrahedron.\marginnote{adjust figure labels. \\ actually, looks ok.}
\begin{figure}[!h]
{\centering

\def\svgwidth{1.1in}
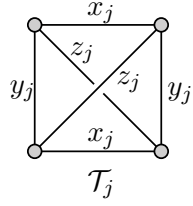

}
\caption{The edges corresponding to the diagonal exchange are labeled $z$.}
\end{figure}

Throughout this section, indices that are out of range will be interpreted cyclically. For example $z_{N+1} := z_1$ and $z_0 := z_{N}$. This convention allows for a much more efficient description of the equations.

\bigskip 

{\noindent \bf Thurston's Equations.}

\smallskip

Many of the edges in the monodromy tetrahedralization meet exactly four faces. This happens when a given edge in $T^2$ lies in two consecutive triangulations $\tau_{j-1}$, $\tau_j$, but does not lie in either $\tau_{j-2}$ or $\tau_{j+1}$. This will be the case if $W_{j} = W_{j-2} R R$, in other words if $\mathcal T_j$ is labeled $RR$.

\noindent In this case, the holonomy around the given edge takes the form 
\begin{equation}
g_j = z_{j-1}z_{j+1}y_j^2  \label{valence-4-hol-R}
\end{equation}
For every $j$ such that $W_{j} = W_{j-2} R R$, the corresponding edge holonomy $g_j$ has the form (\ref{valence-4-hol-R}). Similarly, for every $k$ such that $W_{k} = W_{k-2} L L$, the corresponding edge holonomy has the form
\begin{equation}
g_j = z_{j-1}z_{j+1}x_j^2 = 1 \label{valence-4-hol-L}
\end{equation}
The other edge holonomies can be read off from the hinge tetrahedra. 
A \emph{hinge edge} is an edge $e$ that occurs in more than two consecutive triangulations $\tau_{j-1},\ldots, \tau_{k}$, where we take $p= k-j+2$ to be the maximal number of consecutive $\tau_i$ containing the edge $e$. In this case, $\mathcal T_{j}$ and $\mathcal T_{k}$ are both hinge tetrahedra. Note also that each hinge tetrahedron contains two distinct hinge edges. The edge $e$ is common to the tetrahedra $\mathcal T_{j-1}, \mathcal T_{j}, \ldots, \mathcal T_{k}, \mathcal T_{k+1}$. In $\mathcal T_{j-1}$, $e$ corresponds to the top edge of the diagonal exchange. In $\mathcal T_{k+1}$, $e$ corresponds to the bottom edge of the diagonal exchange. In the case $\mathcal T_j$ is an $LR$ hinge, we have $$W_k = W_{j-2}L RRR\ldots R L$$ and the edge holonomy for $e$ is given by 
\begin{equation}
g_j = z_{j-1} x_{j}^2 x_{j+1}^2\ldots x_k^2 z_{k+1}.\label{hinge-hol-x}
\end{equation}

\noindent If $\mathcal T_j$ is an $RL$ hinge, then we have $$W_k = W_{j-2}R LLL\ldots L R$$ and the edge holonomy for $e$ is given by 
\begin{equation}
g _j = z_{j-1} y_{j}^2 y_{j+1}^2\ldots y_{k}^2 z_{j+1} = 1. \label{hinge-hol-y}
\end{equation}
Every edge in the monodromy tetrahedralization has an edge holonomy expression which is either of the form (\ref{valence-4-hol-R}), (\ref{valence-4-hol-L}) if the edge is valence four or of the form (\ref{hinge-hol-x}), (\ref{hinge-hol-y}) if the edge is hinge.

\begin{figure}[p]
\begin{center}
    \def\svgwidth{4.2in}
	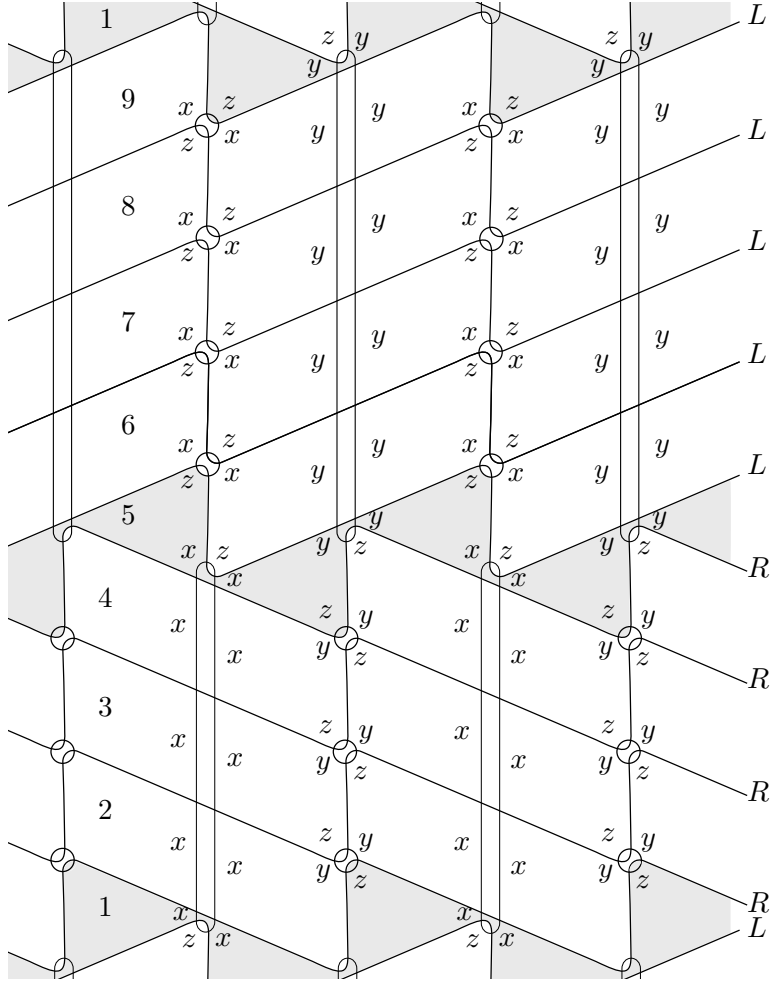
  \end{center}
  \caption[The tessellation of the torus at infinity.]{The edge consistency equations can be read off from a picture of the induced triangulation of $\partial M$. This figure, drawn in the style of Segerman \cite{Segerman-11}, depicts the case $W = R^4 L^5$. The circles and long ovals each represent a vertex of the triangulation. One should imagine the long ovals collapsed down to a point, so that the adjacent quadrilaterals become a fan of triangles around the vertex. The picture is four-periodic going left to right. At any given level of this diagram, the four triangles that touch all come from the same tetrahedron. \label{torus-infinity}}
\end{figure}

All ideal vertices of the tetrahedra $\mathcal T_k$ are identified with one another. The link of the ideal vertex gives a triangulation of $\partial M$. The edge consistency equations can be read off directly from a picture of this triangulation. Vertices in $\partial M$ correspond to edges in $M$. The interior angles of the triangles in $\partial M$ are labeled with the shape parameters corresponding to the edges of the associated tetrahedra in $M$. Figure~\ref{torus-infinity} gives a picture of the combinatorics of $\partial M$ in the case that $W = R^4L^5$.

To conclude this section we summarize the edge consistency equations as follows:
\begin{Proposition}
\label{eqns-punctured-torus}
Let $\phi: T^2 \rightarrow T^2$ be an Anosov map which is conjugate to $$W = R^{m_1}L^{n_1}R^{m_2}L^{n_2}\cdots R^{m_k}L^{m_k}.$$
Then Thurston's edge consistency equations for the canonical ideal triangulation of $M_\phi$ associated to $W$ are described as follows\marginnote{formatting: make sure figure does not split this proposition in half}:

Thinking of $W$ as a string of $R$'s and $L$'s, let $\{j,\ldots,k = j+m_p-1\}$ be the indices of a maximal string of $m_p$ $R$'s. The corresponding $m_p$ equations are:
\begin{flalign*}
& 1 = g_j := z_{j-1} x_j^2 x_{j+1}^2 \cdots x_{k+1}^2 z_{k+2} & \text{(R-fan)}\\
\text{and for each } q = j+1,\ldots,k \ \ \ \ \ \ \ \ \ \ \ \ \ & 1 = g_q := z_{q-1}y_q^2 z_{q+1}    & \text{(R-4-valent)}.
\end{flalign*}
Let $\{j,\ldots,k = j+n_p-1\}$ be the indices of a maximal string of $n_p$ $L$'s. The corresponding $n_p$ equations are:
\begin{flalign*}
& 1 = g_j := z_{j-1} y_j^2 y_{j+1}^2 \cdots y_{k+1}^2 z_{k+2}  & \text{(L-fan)}\\
\text{ and for each } q = j+1,\ldots,k  \ \ \ \ \ \ \ \ \ \ \ \ \ & 1 = g_q := z_{q-1}x_q^2 z_{q+1}  & \text{(L-4-valent)}.
\end{flalign*}
\end{Proposition}
For notation purposes, we write these equations in terms of $\{x_j,y_j,z_j\}$. However, we remind the reader that for each $j = 1,\ldots, N$, $x_j = \frac{z_j-1}{z_j}$, and $y_j = \frac{1}{1-z_j}$. Thus, we think of these equations as depending on the $N$ variables $\{z_j\}$.

\subsection{The real deformation variety}

We look for solutions to the equations of Proposition~\ref{eqns-punctured-torus} over $\mathbb R$ that represent transversely hyperbolic foliations. Requiring that the total dihedral angle around each edge be $2\pi$ amounts to requiring that exactly two of the shape parameters appearing in each equation be negative (see Section~\ref{flattened-tets}). Recall that in the equations of Proposition~\ref{eqns-punctured-torus}, $x_j = \frac{z_j-1}{z_j}$ and $y_j = \frac{1}{1-z_j}$, so that in particular $x_jy_jz_j = -1$ and exactly one of $x_j,y_j,z_j$ lies in each of the components of $\mathbb R \setminus \{0,1\}$. A real shape parameter which is negative is said to have dihedral angle $\pi$, while a positive shape parameter is said to have dihedral angle $0$.

The construction of the monodromy triangulation involved stacking tetrahedra in $T^2 \times \mathbb R$, with each tetrahedron corresponding to a diagonal exchange. From this picture, it would be natural to guess that the edges with dihedral angle $\pi$ should be the edges corresponding to diagonal exchanges, which are labeled $z_j$. This, however, is not the case.

\begin{Proposition}
There is no solution to the equations of Proposition~\ref{eqns-punctured-torus} with all $z_j < 0$.
\end{Proposition}
\begin{proof}
Suppose all $z_j < 0$. Then for all $j$, $x_j > 1$ and $0 < y_j < 1$. Take all the equations involving $x_j$'s and multiply them together. The result is the following:
\begin{equation*}
\prod_{j=1}^N x_j^2 \cdot \prod_{j=1}^N z_j^{\epsilon_j} = 1
\end{equation*}
where each $\epsilon_j = 0,1,$ or $2$. This implies that
\begin{equation*}
\prod_{j=1}^N x_j^{2-\epsilon_j} = \pm\prod_{j=1}^N y_j^{\epsilon_j}
\end{equation*}
which is a contradiction, since the left hand side must be greater than one, while the right hand side must be less than one.
\end{proof}

Due to the structure of the equations, having one of the $z_{j}$ positive actually implies that many other $z_j$'s will be positive as well. In many cases, it can be shown that \emph{all} $z_j$ must be positive. Therefore, a natural subset of solutions to consider is:
\begin{equation}
\Vpl = \{ \text{ real solutions to the equations of Proposition~\ref{eqns-punctured-torus} with } z_j > 0 \text{ for all $j$ }\}
\end{equation}
This set is a union of connected components of the deformation variety. It is also a semi-algebraic set. There are only two possible assignments of dihedral angles (signs) for $\Vpl$:

\begin{Proposition}
Consider an element $(z_1, \ldots, z_N)$ of $\Vpl$. Then $y_j < 0$ if $\mathcal T_j$ is $RR$, and $x_k < 0$ if $\mathcal T_k$ is $LL$. If $\mathcal T_j$ is a hinge tetrahedron, then one of the following two cases holds:
\begin{enumerate}
\item \label{x1-negative} $x_j < 0$ if $\mathcal T_j$ is an $LR$ hinge tetrahedron. $y_k < 0$ if $\mathcal T_k$ is an $RL$ hinge tetrahedron.
\item \label{y1-negative} $x_j < 0$ if $\mathcal T_j$ is an $RL$ hinge tetrahedron. $y_k < 0$ if $\mathcal T_k$ is an $LR$ hinge tetrahedron.
\end{enumerate}
\label{V+-2cases}
\end{Proposition}

\begin{proof}
Begin with the tetrahedron $\mathcal T_1$ which is an $LR$ hinge tetrahedron. Since $z_1 > 0$, we must have $x_1 < 0 $ or $y_1 < 0$. Assume, as in case~\ref{x1-negative}, that $x_1 < 0$. Since $x_1$ appears twice in the first fan equation, choosing $x_1 < 0$ forces all other terms in the first R-fan equation, $$z_N x_1^2 x_2^2 \cdots x_{m_1+1}^2z_{m_1+2} = 1,$$ to be positive (by the $2\pi$ total dihedral angle condition). In particular, $x_{m_1+1} > 0$. Thus, as $z_{m_1+1} > 0$, we must have $y_{m_1 + 1} < 0$. Note that $\mathcal T_{m_1 +1}$ is the second hinge tetrahedron, of type $RL$. Examining the second fan equation (this one is an L-fan), $$z_{m_1} y_{m_1+1}^2 \cdots y_{m_1+m_2+1}^2z_{m_1+m_2+2} = 1,$$ we find that since $y_{m_1+1} < 0$, we have in particular that $y_{m_1+m_2+1} > 0$ and therefore $x_{m_1+m_2+1} < 0$. Note that $\mathcal T_{m_1 + m_2 +1}$ is the third hinge tetrahedron, of type $LR$. This process continues to determine the sign of all hinge shape parameters to be as in case~\ref{x1-negative}. It then follows that the signs of all shape parameters for $RR$ and $LL$ tetrahedra are determined as specified in the Proposition as well.

Similarly, if we begin by choosing $y_1 < 0$, the signs of all other shape parameters are determined to be as in case~\ref{y1-negative}.
\end{proof}

We will focus on the behavior of $\Vpl$. Over the course of the next four sub-sections, we show the following 
\begin{enumerate}
\item $\Vpl$ is smooth of dimension one.
\item All tangent vectors to $\Vpl$ have positive entries (or negative entries).
\item $\Vpl$ is non-empty. In particular $\Vpl$ contains a canonical solution coming from the Sol geometry of the torus bundle obtained by Dehn filling the puncture curve in $M$.
\item $\Vpl$ is locally parameterized by the exponential length of the puncture curve.
\end{enumerate}

The particularly nice form of the equations allows us to prove these properties with relatively un-sophisticated methods.

\subsection{$\Vpl$ is smooth of dimension one.}
We assume case \ref{x1-negative} of Proposition~\ref{V+-2cases}. The other case is symmetric. So, we have 
\begin{itemize}
\item $x_j < 0$ if and only if either $\mathcal T_j$ is an $LR$ hinge tetrahedron or $\mathcal T_j$ is $LL$. 
\item $y_k < 0$ if and only if either $\mathcal T_k$ is an $RL$ hinge tetrahedron or $\mathcal T_k$ is $RR$.
\end{itemize} 
Recall that the edge holonomy expressions, described in Proposition~\ref{eqns-punctured-torus}, are enumerated according to the index of the first $x_j^2$ or $y_j^2$ term appearing in the equation\marginnote{remove repeat equations?}:
\begin{eqnarray*}
g_1(z_1,\ldots,z_N) &:=& z_Nx_1^2x_2^2\ldots x_{m_1+1}^2z_{m_1+2}\\
g_2(z_1,\ldots,z_N) &:=& z_1y_2^2z_3\\
&\vdots&\\g_{m_1}(z_1,\ldots,z_N) &:=& z_{m_1}y_{m_1+1}^2y_{m_1+2}^2\ldots y_{m_1+m_2+1}^2 z_{m_1+m_2+2}\\
g_{m_1+1}(z_1,\ldots,z_N) &:=& z_{m_1}x_{m_1+1}^2z_{m_1+2}\\ &\vdots&
\end{eqnarray*}
The edge consistency equations are given by $g_j = 1$ for all $j = 1, \ldots, N$. In order to determine smoothness and the local dimension, we work with the differentials $dg_j$ of these expressions. Actually, it will be more convenient to work with $ d\log g_j = d g_j / g_j$.  
For convenience we note the differential relationship between $x,y,z$ (leaving off the indices):  \begin{align*} d\log z &= -\frac{1}{y} d\log x = -x d\log y, & d\log y &= -z d\log x. \end{align*}
We choose the following convenient basis for the cotangent space $\mathbb R^{N*}$ at our point $(z_1, \ldots, z_N) \in \Vpl$.  
For indices $j$ such that $x_j < 0$, define $\xi_j = d\log x_j$, $c_j = z_j$, and $t_j = (1-z_j)$ so that  \begin{align*} d\log x_j &= \xi_j,& d\log y_j &= -c_j \xi_j,& d\log z_j &= -t_j \xi_j.\end{align*}
For indices $j$ such that $y_j < 0$, define $\xi_j = d\log y_j$, $c_j = 1/z_j$, and $t_j = (1-z_j)$ so that \begin{align*} d\log x_j &= -c_j\xi_j,& d\log y_j &= \xi_j,& d\log z_j &= -t_j \xi_j.\end{align*}
Note that in both cases $0 < c_j,t_j < 1$. The differential of a fan equation is given by
\begin{eqnarray}
d\log g_j &=& -t_{j-1}\xi_{j-1} + 2\xi_j -2(c_{j+1}\xi_{j+1} + \ldots + c_{k+1}\xi_{k+1}) - t_{k+2}\xi_{k+2}
\label{dlogg1}
\end{eqnarray}
while the differential of a 4-valent equation is given by
\begin{eqnarray}
d\log g_q &=& -t_{q-1} \xi_{q-1} + 2\xi_{q} -t_{q+1}\xi_{q+1}.
\label{dlogg2}
\end{eqnarray}

The kernel of the map $(d\log g_1, \ldots, d\log g_N): \mathbb R^N \rightarrow \mathbb R^N$ is the Zariski tangent space to $\Vpl$. Using the dual basis to $\{ \xi_j\}$ for the domain we let $A$ be the matrix of this map. The matrix $A$ is nearly block diagonal, having a block for each string of $R$'s and a block for each string of $L$'s in the word $W$. A block corresponding to $m_p$ $R$'s will be $m_p \times (m_p+3)$. It overlaps with the following $n_p \times (n_p + 3)$ $L$-block in three columns.

{
\centering

\def\svgwidth{3.0in}

\begingroup
  \makeatletter
  \providecommand\color[2][]{%
    \errmessage{(Inkscape) Color is used for the text in Inkscape, but the package 'color.sty' is not loaded}
    \renewcommand\color[2][]{}%
  }
  \providecommand\transparent[1]{%
    \errmessage{(Inkscape) Transparency is used (non-zero) for the text in Inkscape, but the package 'transparent.sty' is not loaded}
    \renewcommand\transparent[1]{}%
  }
  \providecommand\rotatebox[2]{#2}
  \ifx\svgwidth\undefined
    \setlength{\unitlength}{529.85913086pt}
  \else
    \setlength{\unitlength}{\svgwidth}
  \fi
  \global\let\svgwidth\undefined
  \makeatother
  \begin{picture}(1,0.80134507)%
    \put(0,0){\includegraphics[width=\unitlength]{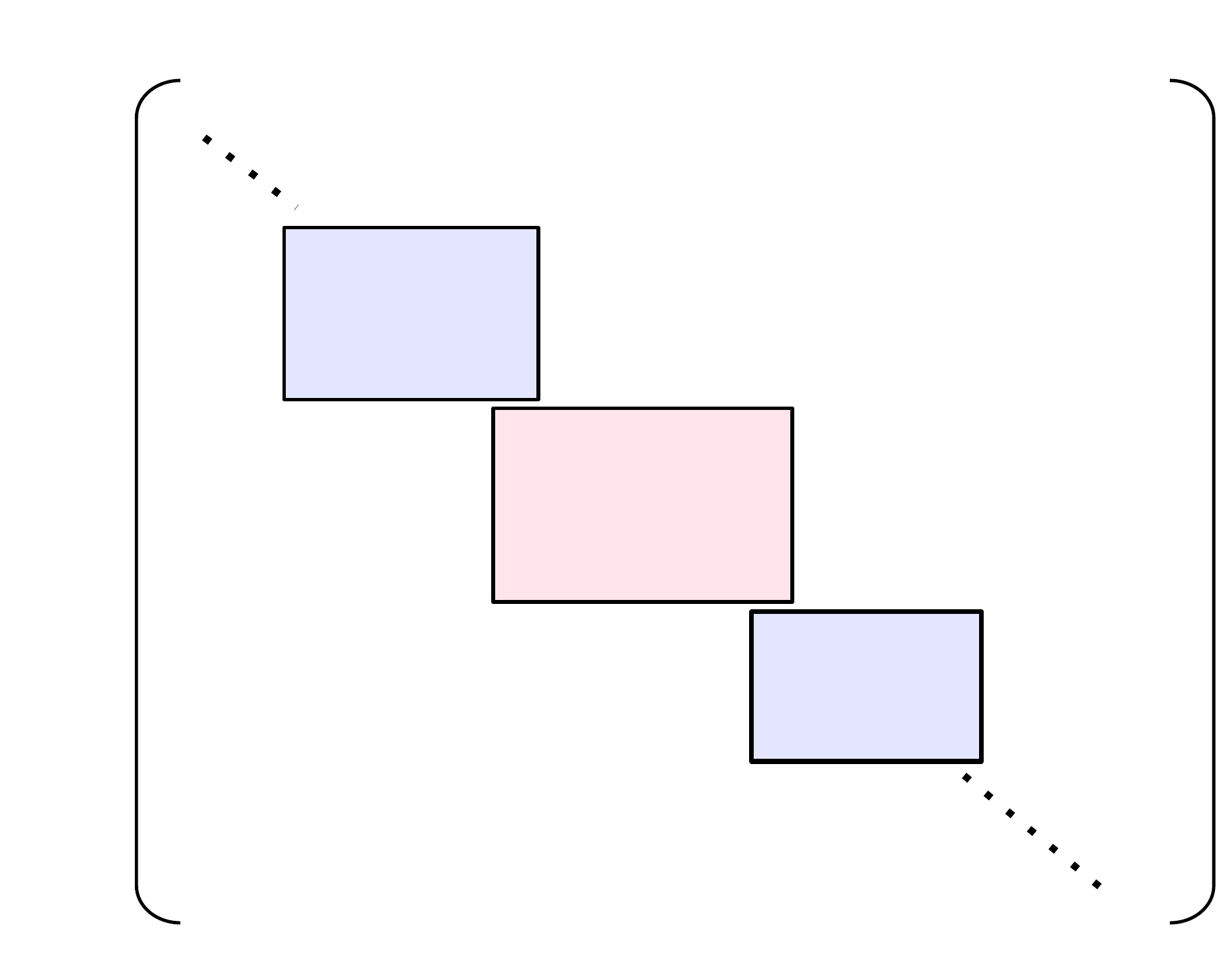}}%
    \put(0.26479115,0.54473961){\color[rgb]{0,0,0}\makebox(0,0)[lb]{\smash{$R$-block}}}%
    \put(0.45718731,0.38254013){\color[rgb]{0,0,0}\makebox(0,0)[lb]{\smash{$L$-block}}}%
    \put(0.63017131,0.23889005){\color[rgb]{0,0,0}\makebox(0,0)[lb]{\smash{$R$-block}}}%
    \put(0.75872299,0.58787778){\color[rgb]{0,0,0}\makebox(0,0)[lb]{\smash{\Huge $0$}}}%
    \put(0.26349703,0.17806537){\color[rgb]{0,0,0}\makebox(0,0)[lb]{\smash{\Huge $0$}}}%
    \put(-0.00137124,0.39850127){\color[rgb]{0,0,0}\makebox(0,0)[lb]{\smash{$A = $}}}%
  \end{picture}%
\endgroup

}

Both $R$-blocks and $L$-blocks have the same form.  Indexing the variables to match Proposition~\ref{eqns-punctured-torus}, each block is described as follows:
\begin{equation*}
\left(\begin{matrix} 
-t_{j-1} & 2 & -2c_{j+1} & -2c_{j+2} & -2c_{j+3} & \ldots & -2c_k & -2c_{k+1} &  -t_{k+2}  \\ 
0 &-t_j & 2 & -t_{j+2} & 0 & \ldots & 0 & 0 & 0 \\
0 & 0 &-t_{j+1} & 2 & -t_{j+3} & \ldots & 0 &0 & 0 \\
0 & 0 & 0 &-t_{j+2} & 2 & \ldots & 0 &0 & 0 \\
& \vdots & & & & \ddots & &\vdots & \\
0 & 0 & 0 & 0 & 0 & \ldots & 2 & -t_{k+1}  & 0\\
\end{matrix}\right)
\end{equation*}
where the $2$'s lie on the diagonal of $A$. For example, if $W = R^4 L^5$ as in Figure~\ref{torus-infinity}, the matrix $A$ is made up of two blocks, an $R$-block of size $3\times 6$ and an $L$-block of size $5 \times 8$ (note: in general the first and last blocks ``spill" over to the other side of the matrix):
\begin{equation*}
A = \begin{pmatrix} 
 2 & -2c_2 & -2c_3 & -2c_4 & -2c_5 & -t_6 & 0 & 0& -t_9  \\ 
-t_1 & 2 & -t_3 & 0 & 0 & 0 & 0 & 0 & 0\\
0 &-t_{2} & 2 & -t_{4} & 0 & 0 &0 & 0 & 0 \\
0 & 0 &-t_{3} & 2 & -t_{5} & 0 & 0 &0 & 0 \\
-2c_1 & -t_2& 0 &-t_{4} & 2 & -2c_6 & -2c_7 &-2c_8 & -2c_9 \\
0 & 0 & 0 & 0 & -t_5 & 2 & -t_7 & 0 & 0 \\
0 & 0 & 0 & 0& 0 &-t_6 & 2 & -t_8 & 0\\
0 & 0 & 0 & 0 &0 &0& -t_7 & 2 & -t_9\\
-t_1 & 0 & 0 & 0 &0 & 0 & 0& -t_8 & 2 \\
\end{pmatrix}
\end{equation*}
We note several important properties of $A$. First, all diagonal entries are equal to $2$. Second, all entries away from the diagonal are non-positive. Third, the entries one off of the diagonal are strictly negative.
Finally, each column sums to zero, which is the differential version of the fact that the product of all $g_j$'s is identically equal to one.

\begin{Proposition}
\label{prop:onedkernel}
The matrix $A$ has one dimensional kernel. \label{one-d-kernel}
\end{Proposition}
\begin{proof}
It will be more convenient to work with the transpose $A^T$, which also has the properties listed above, except
that the rows sum to zero and not the columns.
We write $$A^T = 2I - B-D$$ where $I$ is the $N\times N$  identity matrix, $B = (b_{ij})$ is a matrix with positive entries one off the diagonal and zeros otherwise, and $D = (d_{ij})$\marginnote{ugh, we already have $c$'s above \\ fixed} is a matrix with non-negative entries that is zero within one place of the diagonal. That is $b_{ij}, d_{ij} \geq 0$ for all indices, $b_{ij} > 0$ if and only if $|i-j| = 1$, and $d_{ij} = 0$ if $|i-j| \leq 1$. Now, since the rows of $A^T$ sum to zero, we have that $$v = \begin{pmatrix} 1\\ \vdots \\ 1\end{pmatrix} \ \in \ \ \text{ker} \ A^T.$$
Suppose $w$ is another non-zero vector with $w \in \text{ker} \ A^T$. Then, let $u = w - \text{min}(w)v \in \text{ker }A^T$. Note that all entries of $u$ are non-negative and at least one entry $u_p = 0$. Next, consider the $p^{th}$ entry of $A^T u$:
\begin{eqnarray*}
0 \ = \ -(A^T u)_p &=& -2u_p + (Bu)_p + (Du)_p\\
&=& 0 + b_{p,p-1} u_{p-1} + b_{p,p+1}u_{p+1} + (Du)_p \\
&\geq& b_{p,p-1} u_{p-1} + b_{p,p+1}u_{p+1}.
\end{eqnarray*}
This implies that $u_{p-1} = u_{p+1} = 0$. It follows by induction that $u = 0$ and so $w$ is a multiple of $v$. Thus $A^T$ has one dimensional kernel and so does $A$.
\end{proof}

It follows that $\Vpl$ is smooth and has dimension one. \marginnote{brushing under the rug: it is well-known that the larger complex variety has complex dimension one, at least the irreducible component containing the complete solution does. So $\Vpl$ couldn't have dimension zero (if its part of the same irreducible component...)}

\subsection{Positive tangent vectors}
Actually, we can spiff up the proof of Proposition~\ref{one-d-kernel} to get the following:
\begin{Proposition}
The kernel of $A$ is spanned by a vector with strictly positive entries. 
\label{positive}
\end{Proposition}
\begin{proof}
Again, it is simpler to work with $A^T$.
\begin{Lemma}
The range of $A^T$ does not contain any vectors with all non-negative entries (other than the zero vector).
\end{Lemma}
\begin{proof}
Let $h \in \mathbb R^N$ have non-negative entries and suppose there is $w \in \mathbb R^N$ such that $A^T w = h$. Set $u = w - \min(w)v$, where $v \in \ker A^T$ is, as above, the vector of all $1$'s. Then all entries of $u$ are non-negative, $A^T u = h$, and at least one $u_p = 0$. Following the same argument and notation from the proof of Proposition~\ref{prop:onedkernel} above, we have
\begin{eqnarray*}
0 \ = \ -(A^Tu - h)_p &=&-2u_p + (Bu)_p + (Du)_p + h_p \\
&=& 0 + b_{p,p-1} u_{p-1} + b_{p,p+1}u_{p+1} + (Du)_p + h_p \\
&\geq& b_{p,p-1} u_{p-1} + b_{p,p+1}u_{p+1} + h_p.
\end{eqnarray*}
which shows $u_{p-1}$, $u_{p+1}$, and  $h_p$ are equal to zero. Proceeding inductively, we get that each entry of $u$ and each entry of $h$ is zero.
\end{proof}
The Lemma implies the Proposition as follows. Suppose $u$ spans $\ker A$, and suppose there are two entries $u_i, u_j$ of $u$ such that $u_i u_j \leq 0$. Then, let $v \in \RR^N$ be the vector with entries $v_i = |u_j|$, $v_j = |u_i|$, and all other entries $v_k =0$. Since $v$ is orthogonal to $u$, we have that $v$ is in the range of $A^T$, contradicting the Lemma.
%
\end{proof}

We have (nearly) shown\marginnote{this is just part of Theorem 1, does it need a new name? also, there is one small issue about smoothness: variety could be a point. right?}:
\begin{Proposition}
The tangent space at a point of $\Vpl$ is spanned by a vector with positive components (with respect to $z_j$-coordinates).
\label{prop:pos}
\end{Proposition}
\begin{proof}
Proposition~\ref{positive} gives that the tangent space to the deformation variety is spanned by a vector $u$ in $\mathbb R^N$, whose coordinates with respect to the basis dual to $\{ \xi_j\}$ are positive. Recall that $\xi_j = d\log x_j$ if $x_j < 0$, or $\xi_j = d\log y_j$ if $y_j < 0$. Hence, if $x_j < 0$, then $$dx_j(u) = x_j \xi_j(u)  < 0$$ and if $y_j < 0$, then  $$dy_j(u) = y_j \xi_j(u) < 0.$$ Let $v = -u$.  We remind the reader that 
\begin{align*} dx &= d\left(\frac{z-1}{z}\right) = \frac{1}{z^2} dz& dy = d\left(\frac{1}{1-z}\right) = \frac{1}{(1-z)^2}dz
\end{align*} so that $dz_j(v) > 0$ if and only if $dx_j(v) > 0$ if and only if $dy_j(v) > 0$.
That is, $z_j$ increases in the direction of $v$ if and only if $x_j$ increases in the direction of $v$ if and only $y_j$ increases in the direction of $v$. Hence $v$ has positive coordinates in the standard ($z_j$) basis.
\end{proof}

\subsection{$\Vpl$ is non-empty}
\label{canonical-solutions}
In this section we construct two ``canonical" solutions, $(z_j^+)$ and $(z_j^-)$, to the edge consistency equations which serve as basepoints for the two components of $\Vpl$. These solutions correspond to certain projections of the Sol geometry of the torus bundle associated to $\phi$. We construct them directly from the natural affine $\mathbb R^2$ structure of the layered triangulations used to construct the monodromy tetrahedralization of $M$.

By construction, the punctured torus bundle $M$ comes equipped with a projection map $\pi: \widetilde M \rightarrow \mathbb R^2$ which induces a one dimensional foliation of $M$ with a transverse affine linear structure. Think of $\pi$ as a developing map for the transverse structure. The holonomy $\sigma: \pi_1 M \rightarrow \text{Aff}^+ \mathbb R^2$ can be described as follows: 
\begin{align*}
\sigma(\alpha) &: (x,y) \mapsto (x+1,y),&
\sigma(\beta) &: (x,y) \mapsto (x,y+1),&
\sigma(\gamma) &: (x,y) \mapsto \phi(x,y)
\end{align*}
where $\alpha, \beta$ generate the fiber $\pi_1 T^2$ and $\gamma$ is a lift of the base circle. We can use $\pi$ to project our tetrahedra onto parallelograms in $\mathbb R^2$. Begin by choosing a lift $\tilde{\mathcal T_1}$ of the first tetrahedron. We can choose the lift that projects onto the square $P_1$ with bottom left corner at the origin. We then ``develop" consecutive tetrahedra into $\RR^2$ along a path in $\widetilde M$. The result is a sequence of parallelograms $P_j$, with each consecutive pair overlapping in a triangle. 
\begin{figure}[!h]
{\centering

\def\svgwidth{2.3in}
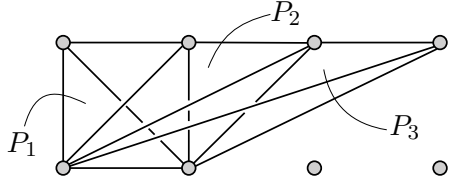

}
\caption{The development of tetrahedra into $\RR^2$ is a union of parallelograms.}
\end{figure}

\noindent The bottom faces of $\tilde{\mathcal T_j}$ map to triangles of $\tilde \tau_{j-1}$ and the top faces map to triangles of $\tilde \tau_j$, where $\tilde \tau_j$ is the lift of the triangulation $\tau_j$ to $\mathbb R^2$. The face glueing maps are realized in $\mathbb R^2$ as combinations of the affine linear transformations $\sigma(\alpha), \sigma(\beta),$ and $\sigma(\gamma)$.

\begin{figure}[!h]
{\centering

\def\svgwidth{3.6in}
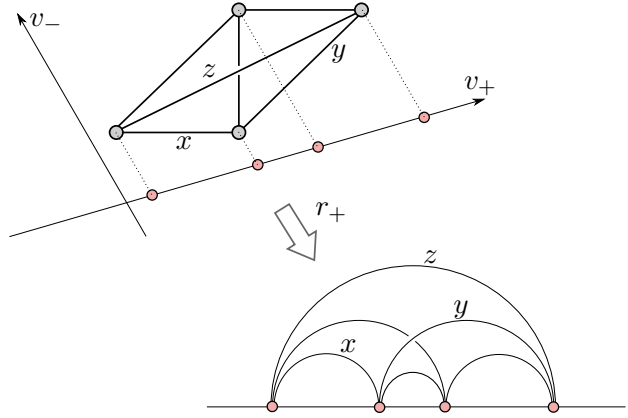

}
\caption{\label{project-tet}The projection $r_+$ determines flattened tetrahedra in $\HH^2$.}
\end{figure}

We now construct $\mathbb H^2$ tetrahedra as follows. Let $v_+, v_-$ be the eigenvectors of $\phi$ corresponding to the eigenvalues $\lambda_+ > 1$, and $\lambda_- < 1$ respectively. Let $r_+, r_-: \mathbb R^2 \rightarrow \mathbb R$ be the coordinate functions with respect to the basis $\{v_+,v_-\}$. For each $j$, project the vertices of $P_j$ to $\mathbb R$ using, e.g., $r_+$. 
The vertices project to four distinct real numbers which we use to build an $\mathbb H^2$ tetrahedron. Orient the resulting $\mathbb H^2$ ideal tetrahedron compatibly with the original tetrahedron. It is an easy exercise to show that this process always produces $\mathbb H^2$ tetrahedra that are folded along the $z$-edges corresponding to diagonal exchanges (i.e. the shape parameter $z$ has zero dihedral angle). See Figure~\ref{project-tet}.

Next, note that $r_+$ takes translations in $\mathbb R^2$ to translations in $\mathbb R$ and $r_+$ converts the action of $\phi$ into scaling by $\lambda_+$ on $\mathbb R$. Hence, the map $r_+\circ \pi : \widetilde M \rightarrow \mathbb R$ is equivariant, converting covering transformations to similarities of $\mathbb R$. In other words, the face glueing maps for our $\mathbb H^2$ tetrahedra are realized by hyperbolic isometries which fix $\infty$. Hence, the shape parameters for these $\mathbb H^2$ tetrahedra are well-defined. 
Using $r_-$ in place of $r_+$ we get a different set of $\HH^2$ shape parameters. The following proposition shows that the condition on dihedral angles is satisfied so that the $r_+$ and $r_-$ shape parameters each determine a solution to Thurston's equation lying in $\Vpl$. It is a corollary of the proof that these solutions lie in distinct components of $\Vpl$.

\begin{Proposition}
\label{r+dihedral}
The collapsed $\HH^2$ triangulations determined by $r_+, r_-$ have total dihedral angle $2\pi$ around each edge in $\mathcal T$.
\end{Proposition}
\begin{proof}
Let $e$ be an edge of the triangulation $\mathcal T$. Recall that $e$ is an edge of consecutive tetrahedra $\mathcal T_{j-1}, \ldots, \mathcal T_{k+1}$, for $k \geq j$ (with $k=j$ if $e$ is $4$-valent, and $k > j$ if $e$ is a hinge edge). In $\widetilde M$, (a lift of) $e$ is bordered by one lift each of $\mathcal T_{j-1}, \mathcal T_{k+1}$, and two lifts of each $\mathcal T_i$ for $i=j,\ldots,k$. The tetrahedra are represented by $2(k-j+2)$ parallelograms $P_i$ which are layered around the corresponding edge $e'$ in $\RR^2$. We number the parallelograms in cyclic order around $e'$ so that $P_1$ is the image of $\mathcal T_{j-1}$, and $P_{k-j+3}$ is the image of $\mathcal T_{k+1}$. For $i=j,\ldots, k$, the two lifts of $\mathcal T_{i}$ that border $e$ map to $P_{i-j+2}$ and $P_{2k-j+4 -i}$, which are translates of one another. Let the endpoints of $e'$ be $p,q \in \ZZ^2$. For each $s = 1,\ldots, 2(k-j+2)$, let $e_s'$ be the edge opposite $e'$ in $P_s$ with endpoints $p_s,q_s$. Note that in the cases $s=1$ and $s=k-j+3$, the edges $e', e_s'$ are diagonals of $P_s$. Let $e_+$ (resp. $e_-$) be the geodesic in $\HH^2$ connecting $r_+(p)$ to $r_+(q)$ (resp. $r_-(p)$ to $r_-(q)$).  For each $s$, the $\HH^2$ ideal tetrahedron $T_s^+$ with vertices $r_+(p), r_+(q), r_+(p_s), r_+(q_s)$ has dihedral angle $\pi$ at $e_s'$\marginnote{changed $e''$ to $e_s'$. should fix in thesis too} if and only if the intervals $r_+(e')$ and $r_+(e_s')$ overlap \emph{partially} (with neither one contained in the other). We will use this characterization to show that the total dihedral angle around $e_+$ is $2\pi$ (and similarly for $e_-$).

The union of the $e_s'$ is a closed polygonal loop in $\RR^2$ with a particularly nice structure. The three edges $e', e_1', e_{k-j+3}'$ share a common midpoint. 
The edges $e_2', \ldots, e_{k-j+2}'$ are each parallel to $e'$; their union forms a straight line segment $I_1 \subset \RR^2$. Similarly, the edges $e_{k-j+4}',\ldots,e_{2(k-j+2)}'$ are each parallel to $e'$ and their union forms a straight line segment $I_2$. Orienting $I_1,I_2$ in the direction of increasing $s$, we have that $I_2$ is a translate of $I_1$ with the same orientation. 
Hence the union of the $e_s'$ is a closed polygonal loop with four straight sides $e_1', I_1, e_{k-j+3}', I_2$. In light of this, the proof will be complete after demonstrating the following lemma.

\begin{figure}[h]
{\centering

\def\svgwidth{3.8in}
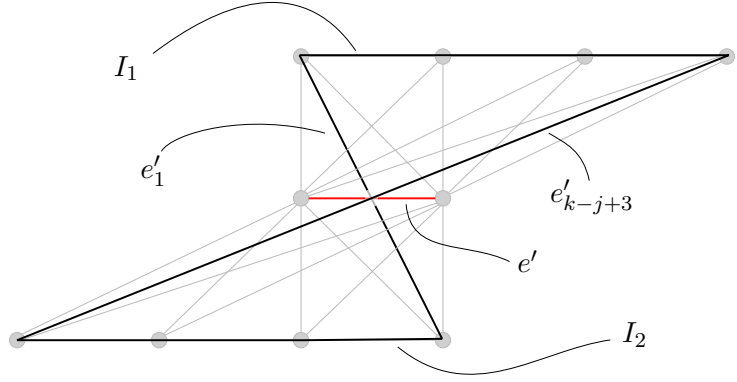

}
\caption{The development of tetrahedra around an edge.}
\end{figure}
\begin{Lemma}
The images of the edges $e_1', e', e_{k-j+3}'$ are nested as follows: 
\begin{align*}
r_+(e_1') &\subset r_+(e') \subset r_+(e_{k-j+3}')\\
r_-(e_{k-j+3}') &\subset r_-(e') \subset r_-(e_1').
\end{align*}
\end{Lemma}
\begin{proof}
Let $P_0$ be the base parallelogram (actually a square) with vertices $(0,0)$, $(1,0)$, $(1,1)$, and $(0,1)$. By construction the parallelogram $P_1$, which corresponds to $\mathcal T_{j-1}$, is given by a translate of $W_{j-2} P_0$ where $W_{s}$ is the product of the first (left-most) $s$ letters in the word $W$ describing $\phi$. The edge $e_1'$ is the bottom diagonal of $P_1$, which is a translate of the vector $W_{j-2}\twovector{1}{-1}$ and the edge $e'$ is the top diagonal of $P_1$, which is a translate of the vector $W_{j-2} \twovector{1}{1}$. Similarly, $P_{k-j+3}$, which corresponds to $\mathcal T_{k+1}$ is a translate of $W_{k} P_0$. So $e'$, which is the bottom diagonal of $P_{k-j+3}$, is a translate of $W_k \twovector{1}{-1}$ and $e_{k-j+3}'$, which is the top diagonal of $P_{k-j+3}$, is a translate of $W_k \twovector{1}{1}$. Recall that either $W_k = W_{j-2}R L^{k-j}R$ or $W_k = W_{j-2} L R^{k-j} L$. From this it is easy to check that $W_{j-2} \twovector{1}{1}= \pm W_k \twovector{1}{-1}$, so they determine the same line segment up to translation in~$\RR^2$. 

Next, we may assume that $v_+$ lies in the positive quadrant and that $v_-$ has negative first coordinate and positive second coordinate (this is easy to check). Hence $W_s \twovector{0}{1}$ and $W_s \twovector{1}{0}$, which lie in the positive quadrant, have $r_+ > 0$. 
Thus we have that for any $s$,
\begin{align*}
r_+ W_s \twovector{1}{1} - r_+ W_s \twovector {1}{-1} &= 2r_+ W_s \twovector{0}{1} &> 0\\
r_+ W_s \twovector{1}{1} + r_+ W_s \twovector {1}{-1} &= 2r_+ W_s \twovector{1}{0} &> 0.
\end{align*}
This implies that $\left| r_+ W_s \twovector{1}{1} \right| > \left| r_+ W_s \twovector{1}{-1} \right|.$
Applying this fact with $s=j-2$ and $s=k$ gives that the lengths of the intervals $r_+(e')$, $r_+(e_1')$, and $r_+(e_{k-j+3}')$ are ordered as follows: 
\begin{align*}
|r_+(e_{k-j+3}')| &> |r_+(e')| > |r_+(e_1')|.
\end{align*}
Thus, as $e'$,$e_1'$, and $e_{k-j+3}'$ share a common midpoint, the $r_+$ statement of the Lemma follows. 

\begin{figure}[!h]
{\centering

\def\svgwidth{4.0in}
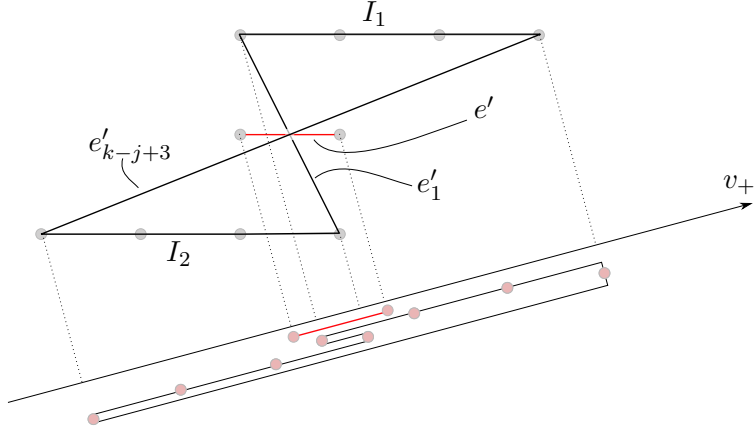

}
\caption{The $r_+$ projection of the edges opposite $e'$.}
\end{figure}

The $r_-$ statement is similar. In this case $r_- W_s \twovector{0}{1} > 0$ while $r_- W_s \twovector{1}{0} < 0$. Thus
\begin{align*}
r_- W_s \twovector{1}{1} - r_- W_s \twovector {1}{-1} &= 2r_- W_s \twovector{0}{1} &> 0\\
r_- W_s \twovector{1}{1} + r_- W_s \twovector {1}{-1} &= 2r_- W_s \twovector{1}{0} &< 0.
\end{align*}
It follows that for any $s$, $\left| r_- W_s \twovector{1}{1} \right| < \left| r_- W_s \twovector{1}{-1} \right|.$
So the lengths of the intervals $r_-(e')$, $r_-(e_1')$, and $r_-(e_{k-j+3}')$ are ordered like so:
\begin{align*}
|r_-(e_{k-j+3}')| &< |r_-(e')| < |r_-(e_1')|.
\end{align*}
The $r_-$ part of the Lemma now follows.

\end{proof}
By the Lemma, and the above characterization of the edges $e_s'$, we must have that $r_+(e_s')$ partially overlaps $r_+(e')$ if and only if $s = 2$ or $s =2(k-j+2).$ Hence the $\HH^2$ tetrahedron $T_s^+$ has dihedral angle $\pi$ at $e_+$ if and only if $s = 2$, or $s = 2(k-j+2)$. Similarly, the $\HH^2$ tetrahedron $T_s^-$ in the $r_-$ tetrahedralization has dihedral angle $\pi$ if and only if $s = k-j+2$ or $s = k-j+4$. Note that this shows that $r_+$ produces a solution in case~\ref{x1-negative} of Proposition~\ref{V+-2cases}, while $r_-$ produces a solution in case~\ref{y1-negative}.
\end{proof}

 We note that $r_+ \circ \pi$ and $r_- \circ \pi$ give maps to $\mathbb R \subset \mathbb{RP}^1$ which are invariant under reducible representations $\rho_+, \rho_- : \pi_1 M \rightarrow \PSL(2,\mathbb R)$ respectively. 
 Let $(z^+_j)$ and $(z^-_j)$ denote the two solutions coming from $r_+$ and $r_-$ respectively. These solutions determine two distinct transversely hyperbolic foliations $\mathcal F^+$ and $\mathcal F^-$ whose holonomy representations are $\rho_+, \rho_-$. These two transversely hyperbolic foliations come from projecting the Sol geometry of $M_\eps$ onto the two vertical hyperbolic planes in Sol, where $M_\eps$ is the torus bundle obtained by Dehn filling $M$ along the puncture curve $\eps$.

\subsection{A local parameter}
\label{subsec:local}
Let $\epsilon$ represent the curve encircling the puncture in $T^2 \subset M$. We show in this section that the length of $\epsilon$ is a local parameter for $\Vpl$. This will follow, after some calculation, from the fact that the tangent direction to $\Vpl$ must increase all shape parameters.

\begin{figure}[h]
{\centering

\def\svgwidth{3.6in}
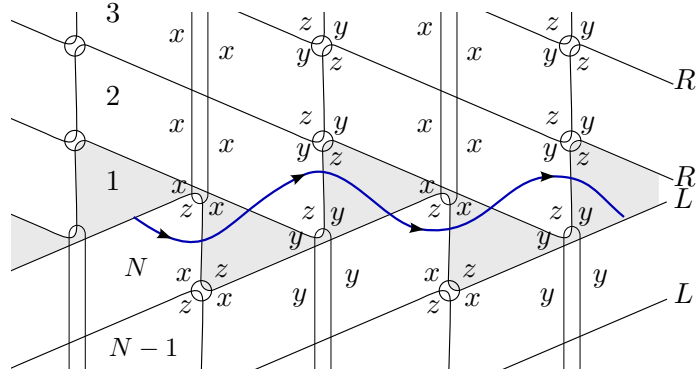

}
\caption{\label{puncture-curve}The puncture curve $\epsilon$ is drawn in blue}
\end{figure}

In order to ease the upcoming computation, we change notation slightly, and write the decomposition of our Anosov map $\phi$ as:
\begin{equation*}
A \phi A^{-1} = W = R^{s_1}L^{s_2}R^{s_3}L^{s_4}\cdots R^{s_{K-1}}L^{s_K}.
\end{equation*}
Let $M_p$ denote the index of the $p^{th}$ hinge tetrahedron, $M_p = 1+\sum_{j=1}^{p-1}s_j,$
where we define $M_{K+1} = N+1$ and $N = \sum s_j$ is the total length of the word $W$. 
We assume (as the other case is similar) that we are in case \ref{x1-negative} of Proposition \ref{V+-2cases}, that is, $x_1 < 0$. We adopt the following notation:
\begin{align*}
\alpha_j &= \left\{ \begin{array}{ll}
         x_j & \mbox{if $x_j < 0$}\\
         y_j & \mbox{if $y_j < 0$},\end{array} \right. &
 \beta_j &= \left\{ \begin{array}{ll}
         x_j & \mbox{if $x_j > 0$},\\
         y_j & \mbox{if $y_j > 0$}.\end{array} \right. 
\end{align*}
Note that, in either case, we have $\alpha_j \beta_j z_j = -1$. 
As usual, the indices $i$ of the $\alpha_i, \beta_i, z_i$ are to be interpreted cyclically so that, for example, $\beta_{N+1} := \beta_1$.
The exponential length of $\epsilon$ can be read off from Figure~\ref{puncture-curve}:
\begin{align*}
H(\epsilon) &= (z_N x_1 z_2^{-1} y_1^{-1} )^2\\
&= (z_N \alpha_1 z_2^{-1} \beta_1^{-1})^2.
\end{align*}
Further, the gluing equations take a nice form with respect to these coordinates. If the $j^{th}$ edge is 4-valent (of type R or L), then the $j^{th}$ edge equation is:  $$1 = g_j =  z_{j-1} \alpha_j^2 z_{j-1}.$$
If $j= M_p$ is the index of the $p^{th}$ hinge edge, then the $j^{th}$ edge equation is: $$1 = g_{M_p} = z_{M_p-1} \alpha_{M_p}^2 \beta_{M_p+1}^2\ldots \beta_{M_{p+1}}^2 z_{1+M_{p+1}}.$$
\begin{Proposition} \label{prop:parameterization}
 $H(\epsilon)$ can be expressed in the following form:
\begin{equation}
H(\epsilon)^{N/2} = \prod_{p=1}^{K} \left( z_{1+M_{p+1}}^{-s_{p}} \prod_{j=1+M_{p}}^{M_{p+1}} \beta_{j}^{-2j+2M_{p}} \right). \label{lots-o-products}
\end{equation}
In particular, $H(\epsilon)$ is a local parameter for $\Vpl$.
\end{Proposition}
\begin{proof}
We consider a product of powers of the gluing equations as follows: 
\begin{eqnarray*}
\prod_{i=1}^{N-1} g_{i+1}^i &=& \prod_{i=1}^{N-1} (z_i\alpha_{i+1}^2z_{i+2})^i \prod_{p = 2}^K \left( z_{1+M_p}^{-(M_p-1)} z_{1+M_{p+1}}^{(M_p-1)}\prod_{j=1+M_p}^{M_{p+1}} \beta_j^{2(M_p-1)}\right)\\
&=& z_1^N z_N^{-N}\prod_{i=1}^{N-1} \left(z_{i+1}^{2i}\alpha_{i+1}^{2i}\right)  \left(\prod_{p = 2}^K z_{1+M_p}^{1-M_p} \prod_{j=1+M_p}^{M_{p+1}} \beta_j^{2M_p-2}\right)\prod_{p=3}^{K+1}\left( z_{1+M_{p}}^{-1+M_{p-1}}\right)\\
 &=& z_1^N z_N^{-N}\prod_{i=1}^{N-1} \beta_{i+1}^{-2i}  \left(\prod_{p = 2}^K\prod_{j=1+M_p}^{M_{p+1}} \beta_j^{2M_p-2}\right)
 z_2^{N}\prod_{p=2}^{K+1} \left( z_{1+M_p}^{M_{p-1}-M_{p}}\right) \\
  &=& z_1^N z_N^{-N}  z_2^{N} \beta_1^{2N}\left(\prod_{p = 1}^K\prod_{j=1+M_p}^{M_{p+1}} \beta_j^{-2(j-1)+2M_p-2}\right)
\prod_{p=2}^{K+1} \left( z_{1+M_p}^{-s_{p-1}}\right) \\
&=& ( z_N^{-1} \alpha_1^{-1}z_2 \beta_1)^{N} \prod_{p=1}^{K} \left( z_{1+M_{p+1}}^{-s_{p}}\prod_{j=1+M_{p}}^{M_{p+1}} \beta_{j}^{-2j+2M_{p}}\right). \\
\end{eqnarray*}
So, as $H(\epsilon) = (z_N \alpha_1 z_2^{-1} \beta_1^{-1})^2$, we have the result.
We have expressed $H(\eps)$ as a product of negative powers of shape parameters. Therefore since all shape parameters increase (or decrease) together along $\Vpl$, we have that $H(\eps)$ is a local parameter.
\end{proof}
\begin{Corollary}\label{cor:parameterize}
The exponential length of the puncture curve, $H(\epsilon)$, can be made arbitrarily small or arbitrarily large .
\end{Corollary}
\begin{proof}
We show that $H(\epsilon)$ can be made arbitrarily small; the second statement is similar. By Proposition~\ref{one-d-kernel}, any point $(z_1,\ldots,z_N) \in \Vpl$ is a smooth point of the deformation variety, and by Proposition~\ref{prop:pos}, all shape parameters can be increased locally. They can be increased globally until some of them go to infinity. Consider such a path. Recall that for some $j$, $\beta_j = x_j$ while for other $k$, $\beta_k = y_k$. We show that one of the $\beta_i$ must approach infinity. Assume not. Then some $z_j \rightarrow \infty$, for recall that all $\alpha_j < 0$. Examining the ${(j+1)}^{st}$ glueing equation (in which $z_j$ appears), we see that $\alpha_{j+1} \rightarrow 0$ as all other terms are positive and increasing. If $\beta_{j+1} > 1$ ($\iff \beta_j = y_j$), then we must have $\beta_{j+1} \rightarrow \infty$. If not then $z_{j+1} > 1$ and $z_{j+1} \rightarrow \infty$. We continue inductively and eventually reach an index $i$ such that $\beta_i > 1$ and $\beta_i \rightarrow \infty$. It is clear from the expression (\ref{lots-o-products}) that $H(\epsilon) \rightarrow 0$.
\end{proof}

\begin{Remark}
It is straightforward to check that the discrete rotational part of the holonomy of $\epsilon$ must be $+2\pi$. 
\end{Remark}

\subsection{The space of transversely hyperbolic foliations}

We conclude the section by describing the image of $\Vpl$ in the space of transversely hyperbolic foliations on~$M$. Two transversely hyperbolic foliations $\mathcal F, \mathcal F'$ are considered equivalent if their submersive developing maps $D, D' : \widetilde M \to \HH^2$ satisfy $D' = g D \widetilde \psi$ where $g \in \operatorname{Isom}\HH^2$ and $\widetilde \psi$ is the lift of a diffeomorphism $\psi$ of $M$.

Recall the two transversely hyperbolic foliations $\mathcal F^+$, $\mathcal F^-$ determined by the solutions $(z^+_j), (z^-_j)$ described in Section~\ref{canonical-solutions}. Let $\Vpl^+$ and $\Vpl^-$ denote the smooth one-dimensional components containing $(z^+_j)$ and $(z^-_j)$ respectively. Huesener--Porti--Su\'arez \cite{Huesener-01} have carefully studied the deformation space of transversely hyperbolic foliations on $M$ near the special points $\mathcal F^+, \mathcal F^-$. Their work implies that, though $\mathcal F^+$ and $\mathcal F^-$ are not equivalent, any structure nearby $\mathcal F^+$ is equivalent to one nearby $\mathcal F^-$. In fact, we argue that $\Vpl^+\setminus \{(z^+_j)\}$ and $\Vpl^-\setminus \{ (z^-_j)\}$ parameterize equivalent transversely hyperbolic foliations. For let $(w^+_j) \in \Vpl^+$ be a solution, different from~$(z_j^+)$, determining a transversely hyperbolic foliation~$\mathcal F$. Denoting the ideal vertices of $(\widetilde M, \widetilde{\mathcal T})$ by $\operatorname{Vert}(\widetilde{\mathcal T})$, then $(w_j^+)$ determines a map $p^+: \operatorname{Vert}(\widetilde{\mathcal T}) \to \RP^1$ equivariant under the holonomy representation of~$\mathcal F$. For each $v \in \operatorname{Vert}(\widetilde{\mathcal T})$, the hyperbolic subgroup representing the stabilizer of~$v$ has two fixed points, $p^+(v)$ and another which we denote~$p^-(v)$. A second solution $(w^-_j) \in \Vpl^-$, giving the same transversely hyperbolic foliation, is determined by the map~$p^-$. Note that for the solution $(z_j^+)$, the holonomy representation is reducible, so that all conjugates of $\pi_1 \partial M$ have a common fixed point $q_0$; in this case $p^-(v) = q_0$ does not correspond to a solution because ideal tetrahedra have four distinct vertices.

\begin{Proposition}
The correspondence $(w^+_j) \mapsto (w^-_j)$ gives a well-defined homeomorphism $\Vpl^+\setminus \{(z^+_j)\} \to \Vpl^-\setminus \{ (z^-_j)\}$.
\end{Proposition}
\begin{proof}
To check well-defined, we must check that $p^-(v) \neq p^-(v')$ whenever $v, v'$ are ideal vertices of the same tetrahedron. Examining the monodromy triangulation, we find that for any two ideal vertices $v, v'$ of the same tetrahedron, there are generators $a , b$ of the punctured torus fiber for which $v$ is fixed by $[a,b]$ and $v'$ is fixed by $a[a,b]a^{-1}$. One easily checks that the holonomy representation corresponding to $(w_j^+)$ maps $[a,b]$ and $a[a,b]a^{-1}$ to hyperbolic elements with distinct fixed points unless the full holonomy representation is reducible; thus the map is well-defined unless $(w_j^+) = (z_j^+)$. Well-defined-ness of the inverse map  $\Vpl^-\setminus \{(z^-_j)\} \to \Vpl^+\setminus \{ (z^+_j)\}$ is similar.
\end{proof}

The proposition shows that the images of $\Vpl^+$ and $\Vpl^-$ in the space of transversely hyperbolic foliations (up to equivalence) on $M$, form the non-Hausdorff space known as the ``line with two origins", the two origins being $\mathcal F^+$ and $\mathcal F^-$. 

\section{$\AdS$ structures with a tachyon}
\label{sec:tachyon}
As in the previous section, let $M$ be a punctured torus bundle with Anosov monodromy equipped with the monodromy triangulation $\mathcal T$. We now apply the results of Section~\ref{subsec:local} to prove Theorem~\ref{thm:tachyons} from the introduction. 
We find $\AdS$ structures on $(M,\mathcal T)$ so that the geometry extends, with a tachyon singularity, to the manifold $M_\eps$ obtained by Dehn filling the puncture curve $\eps$. By the arguments of Section~\ref{AdS-tetrahedra} (see Example~\ref{fig8-AdS}), such a structure is produced by finding a positively oriented, space-like solution to Thurston's equations over $\Rtau$ with the added condition that the holonomy around the puncture curve $\eps$, has exponential ($\Rtau$)-length given by:
\begin{align} \label{eqn:8-tachyon}
H(\eps) &= e^{\tau \phi}
\end{align}
where $\phi$ is the tachyon mass. 
Let $(\lambda_j), (\mu_j)$ be two real solutions to Thurston's equations lying in the same component of $\Vpl$ (say $\Vpl^+$). By Proposition~\ref{AdS-RcrossR}, $(\lambda_j), (\mu_j)$ determine a positive oriented solution over $\Rtau$ if and only if $\lambda_j > \mu_j$ for all $j = 1, \ldots, N$.  
As in Example~\ref{fig8-AdS}, Equation~(\ref{eqn:8-tachyon}) is equivalent to 
\begin{align*}
H_\lambda(\eps) &= e^\phi & H_\mu(\eps) = e^{-\phi}
\end{align*}
where $H_\lambda$ and $H_\mu$ refer to the real exponential length functions for the solutions $(\lambda_j)$ and $(\mu_j)$ respectively. By Proposition~\ref{prop:parameterization}, $\phi$ parameterizes such pairs of solutions $(\lambda_j), (\mu_j)$. By the proof of Proposition~\ref{prop:parameterization}, we have that $\lambda_j > \mu_j$ for all $j$ if and only if the tachyon mass $\phi < 0$. By Corollary~\ref{cor:parameterize}, the tachyon mass $\phi$ can take any value in $(-\infty, 0)$. This completes the proof of Theorem~\ref{thm:tachyons}. 
We note that choosing $(\lambda_j), (\mu_j)$ in the other component of $\Vpl$ will produce the same family of $\AdS$ structures, but with the triangulation spinning in the opposite direction around the singular locus (see discussion at the end of Section~\ref{canonical-solutions}). 

\begin{figure}[h]
{\centering

\def\svgwidth{4.5in}
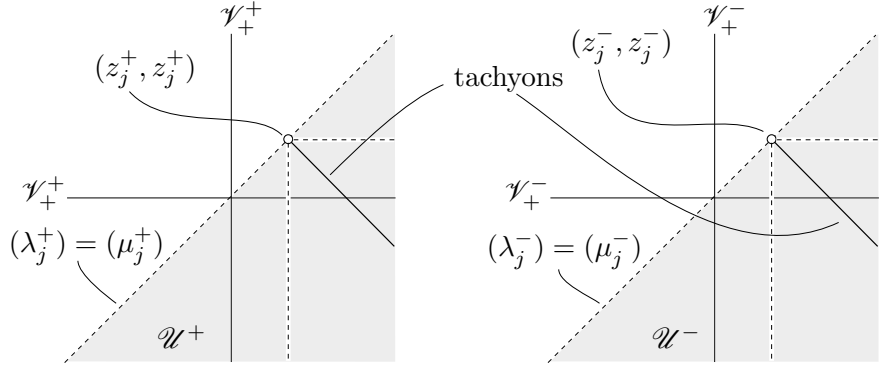

}
\caption{The deformation space of $\AdS$ structures on $M$. The half-spaces $\mathscr U^+$ and $\mathscr U^-$ are identified along the complement (shaded) of the vertical and horizontal dotted lines. The result is not Hausdorff.}
\label{fig:non-Hausdorff}
\end{figure}

Finally we remark that the space of $\AdS$ structures on $M$, without restriction on the geometry at the boundary, inherits the non-Hausdorff behavior of the space of transversely hyperbolic foliations. For, the two sets
\begin{align*}
\mathscr U^+ &= \left\{ (\lambda_j^+,\mu_j^+)_{j=1}^N : (\lambda_j^+), (\mu_j^+) \in \Vpl^+, \text{ and } \lambda_j^+ > \mu_j^+ \text{ for all j} \right\}\\
\mathscr U^- &= \left\{ (\lambda_j^-,\mu_j^-)_{j=1}^N : (\lambda_j^-), (\mu_j^-) \in \Vpl^-, \text{ and } \lambda_j^- > \mu_j^- \text{ for all j} \right\}
\end{align*}
determine smooth two-dimensional subspaces, each embedded in the deformation space of $\AdS$ structures on $M$.
However, a structure determined by $(\lambda_j^+, \mu_j^+) \in \mathscr U^+$ is isomorphic to a structure determined by $(\lambda_j^-, \mu_j^-)\in \mathscr U^-$ (and vice versa) provided that $(\lambda_j^+)$ and $(\mu_j^+)$ are both different from the solution $(z_j^+)$, and $(\lambda_j^-)$ and $(\mu_j^-)$ are both different from the solution $(z_j^-)$.  Thus the non-Hausdorff space obtained by idenitfying $\mathscr U^+$ and $\mathscr U^-$ along the complement of two lines is embedded in the space of $\AdS$ structures on $M$. See Figure~\ref{fig:non-Hausdorff}.

\bibliographystyle{amsalpha}
\small
\bibliography{/Users/jeffdanciger/Desktop/Desktop/bib/mybib}

\end{document}